\documentclass [11pt,oneside]{article}
\usepackage{a4,latexsym,exscale,theorem,epsfig}
\usepackage{amssymb,psfrag,epsf,amsmath,verbatim,bbm,float}
\usepackage{verbatim}
\usepackage{todonotes}
\usepackage{listings}
\usepackage{color}
\usepackage{mathbbol}
\newtheorem{lemma}{Lemma}[section]

\newtheorem{remark}[lemma]{Remark}
\newtheorem{example}[lemma]{Example}

\renewcommand{\theequation}{\arabic{section}.\arabic{equation}}
\renewcommand{\thefigure}{\arabic{section}.\arabic{figure}}

\begin{document}
\newenvironment{proof}{\hspace{-\parindent}{\bf{}Proof:}}
{\hspace*{\fill}$\blacksquare$}
\newcommand {\eps} {\varepsilon}
\newcommand {\Z} {\mathbb{Z}}
\newcommand {\R} {\mathbb{R}}
\newcommand {\N} {\mathbb{N}}
\newcommand {\ang} {\measuredangle}
\newcommand {\e} {{\rm{e}}}
\newcommand {\rank} {{\rm{rank}}}
\newcommand {\Span} {{\rm{span}}}
\newcommand {\card} {{\rm{card}}}
\newcommand {\OO} {\mathcal{O}}
\newcommand {\nr} {\mathcal{N}}
\newcommand{\cZ} {\mathcal{Z}}
\newcommand {\Sn}[1] {\mathcal{S}^{#1}}
\newcommand {\range} {\mathcal{R}}
\newcommand {\kernel} {\mathcal{N}}
\newcommand{\one}{\mathbb{1}}
\renewcommand{\thefootnote}{\fnsymbol{footnote}}

\newcommand{\sect}%
{
  \setcounter{equation}{0}%
 \setcounter{figure}{0}%

  \section
  }



\title{The Right Space for Dynamics: \\
Numerics with Diffeomorphism Equivariance}

\author{Wolf-J\"urgen Beyn\footnotemark[1]\quad
Michail E. Kavousanakis \footnotemark[2] \quad
Yannis G. Kevrekidis \footnotemark[3] 
}
\footnotetext[1]{Department of Mathematics, Bielefeld University 
\texttt{beyn@math.uni-bielefeld.de}}
\footnotetext[2]{School of Chemical Engineering, National Technical University, 
Zografou, Athens, 
\qquad \texttt{mihkavus@chemeng.ntua.gr}}
\footnotetext[3]{Department of Applied Mathematics and
Statistics, Johns Hopkins University, \qquad \texttt{yannis@jhu.edu}}

\maketitle


\begin{abstract}

Among many (equivalent, via invertible transformations) representations of the evolution of a dynamical system, which one is to be preferred? Here we show how the use of infinite-dimensional Lie group theory (and its numerical implementation) allows us to single out one representation, by selecting an element of the group of diffeomorphisms acting on the dynamical system. We present and discuss several types of ``phase conditions" defining the selected representation, and illustrate their computational implementation. Study of dynamics {\em modulo} diffeomorphisms ``liberates" mathematical modeling of physical phenomena from a user's preferred coordinates, and spontaneously selects a ``right latent space" for the system.
\end{abstract}

\noindent \textbf{Keywords:}
Dynamical systems, numerics of PDEs, diffeomorphism equivariance, phase conditions

\noindent \textbf{MSC codes:} 37C05, 37C80, 58D05, 58D19, 65M99

\sect{Introduction}
This paper deals with the solution of evolution equations of the form
\begin{equation} \label{e1}
\frac{d}{dt}u(t)= F(u(t)), \;t\ge 0, \quad u(0)=u_0,
\end{equation}
where the state $u(t)$ at time $t$ and the initial data $u_0$ are
$C^{\infty}$-maps from a finite dimensional manifold $M$ into some finite, $d$-dimensional  space,\ i.e.\ $u(t),u_0 \in C^{\infty}(M,\R^d)$. The manifold $M$ is the coordinate space, whose space carries the independent variables $x$.
Whenever we speak of spatial derivatives we mean derivatives with respect to these coordinates on $M$.

We assume that the system is equivariant with respect to
transformations in coordinate space,~ i.e.
 the vector field $F$ on the right-hand side of \eqref{e1} commutes with the  action $a$ of the group of $C^{\infty}$-diffeomorphisms
$\mathrm{Diff}(M)$ defined by
\begin{equation*} \label{e1.2}
  a(g,u) = u\circ g^{-1}, \quad g \in \mathrm{Diff}(M), \quad u \in C^{\infty}(M,\R^d).
\end{equation*}
We write the equivariance explicitly as 
\begin{equation*} \label{eq1:equivariance}
  F(a(g,u))=a(g,F(u)), \quad g\in \mathrm{Diff}(M), \quad u \in C^{\infty}(M,\R^d)
\end{equation*}
and refer to Remark \ref{rem2} for a more general form.

Symmetries of this kind are abundant in mathematical physics and appear in
quite a variety of field theories (see e.g.\ \cite{MEF72}, \cite[Ch.9]{MR99} and the review \cite{Sch10}).
The main reason for diffeomorphism equivariance of such theories is to formulate physical laws in
a coordinate free form. Certainly, the  most prominent example is 
Einstein's field equation of general relativity \cite{FM72}, \cite{Baumgarte10},  \cite{CB15}
which is written in spacetime geometry - and is thus more general than
the evolution form \eqref{e1}.
Further examples of diffeomorphism equivariant equations \eqref{e1}
are mean curvature flow (\cite{Br78},\cite{Sinestrari2010}), the curve shortening equation (\cite{AAAW13}), the vortex filament equation (\cite{Kh13}),
Maxwell's equation from electrodynamics and the Maxwell-Vlasov equation
from plasma physics (see \cite{MEF72}, \cite[Ch.9]{MR99}, \cite{Sch10}).
For the latter field equations, note that 
equivariance only holds  with respect to {\em an infinite dimensional
subgroup of diffeomorphisms} (gauge and canonical transformations).
Suitable transformations of the underlying domain have also been used in numerical
analysis to design so-called adaptive moving mesh methods; see \cite{huru11}.

Symmetries underly a large variety of methods for solving evolution
equations, both theoretically and numerically. They suggest special coordinate
systems or suitable phase spaces and lead to reduction procedures 
and conserved quantities (Noether invariants); see \cite{FM72}, \cite[Ch.4]{Baumgarte10}.
They also play a prominent role in the bifurcation analysis of steady state
equations  \cite{GS85}, to mention just a few applications.

Our goal in this paper is to develop a numerical approach
which takes advantage of diffeomorphism equivariance by first reformulating \eqref{e1} in
an extended phase space and then discretizing afterwards. More precisely, we rewrite equation
\eqref{e1} for a new state variable $v(t) \in C^{\infty}(M,\R^d)$ with
\begin{equation*}\label{e1a}
  v(t) = a(g(t),u(t))= u(t)\circ g^{-1}(t), \quad g(t)\in \mathrm{Diff}(M),
\end{equation*}
where $g(t)$ is a time-dependent and a-priori unknown diffeomorphism.
In section \ref{sec1} we show that this is equivalent to solving an extended evolutionary
system
\begin{align} \label{e1b}
  \frac{d}{dt}v(t)= & F(v(t))+ \partial_xv(t) \mu(t) ,\quad v(0)=u_0, \\
  \label{e1c}
\frac{d}{dt} g(t) = & \partial_xg(t) \mu(t,\cdot), \quad g(0)= \one. 
\end{align}
for $v(t)\in C^{\infty}(M,\R^d)$, $g(t)\in \mathrm{Diff}(M)$ along with still another time-dependent and unknown vector field
$\mu(t)\in C^{\infty}(M,TM)$.
To fix these extra unknowns and thus to close the system \eqref{e1b}, \eqref{e1c},
we will derive a set of side constraints called {\em phase conditions}. These
are based on minimization principles designed to facilitate the (numerical) solution
of the system. In Section \ref{sec2} we discuss two types of functionals which are
either minimized with respect to the group variable $g$ or with respect to
the vector field $\mu$. In both cases the phase conditions turn out  to be an infinite
set of pointwise constraints. Altogether, we arrive at a {\em partial differential-algebraic
equation} (PDAE), the solution of which
is expected to give better insight into the structure and/or into the longtime dynamics of the
original system \eqref{e1}.

On an abstract level,
the approach generalizes the well-known concept of 
symmetry reduction which appears under various names and in 
various realizations in
the literature. We mention the concept of reduction and reconstruction  by  Rowley et.al.~\cite{RKML03}, the Freezing
method by Beyn et.al.~\cite{BT04,BOR14} and the
method of slices by Froehlich, Cvitanovi\'{c} \cite{FC12}.
The names reflect different aspects of the method,
and for brevity in this paper we adopt the use of
the single word ``Freezing" (but please note that our primary concern is
{\em not} to make dynamic features as stationary as possible, as the word colloquially suggests).

The main goal of this paper is to go one step further from finite to {\em infinite dimensional Lie groups}.
Since arbitrary coordinate systems are involved, this
step turns out to be a major one, both theoretically and numerically.

\begin{figure}
    \begin{tabular}{cc}
    \includegraphics[width=0.5\linewidth]{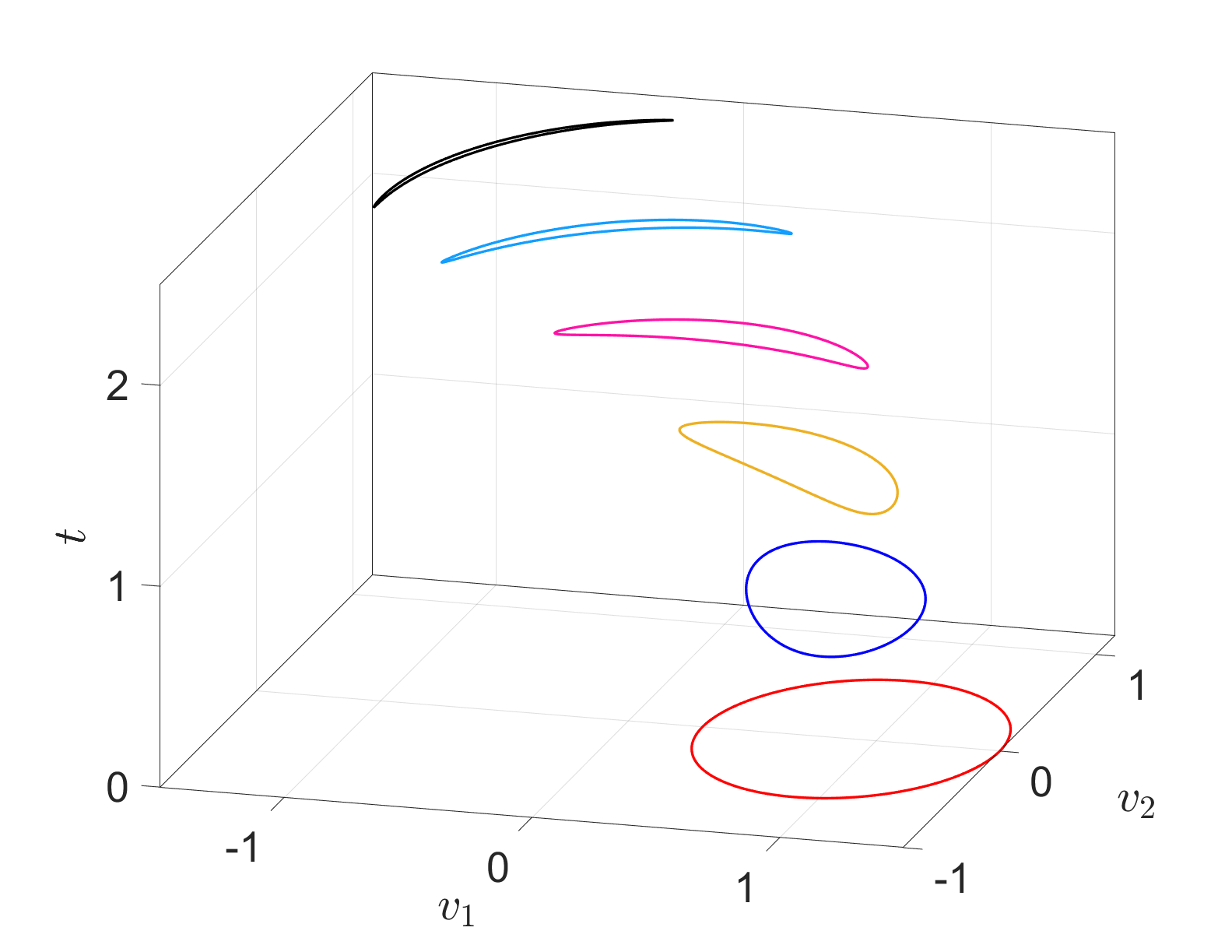}     &  \includegraphics[width=0.5\linewidth]{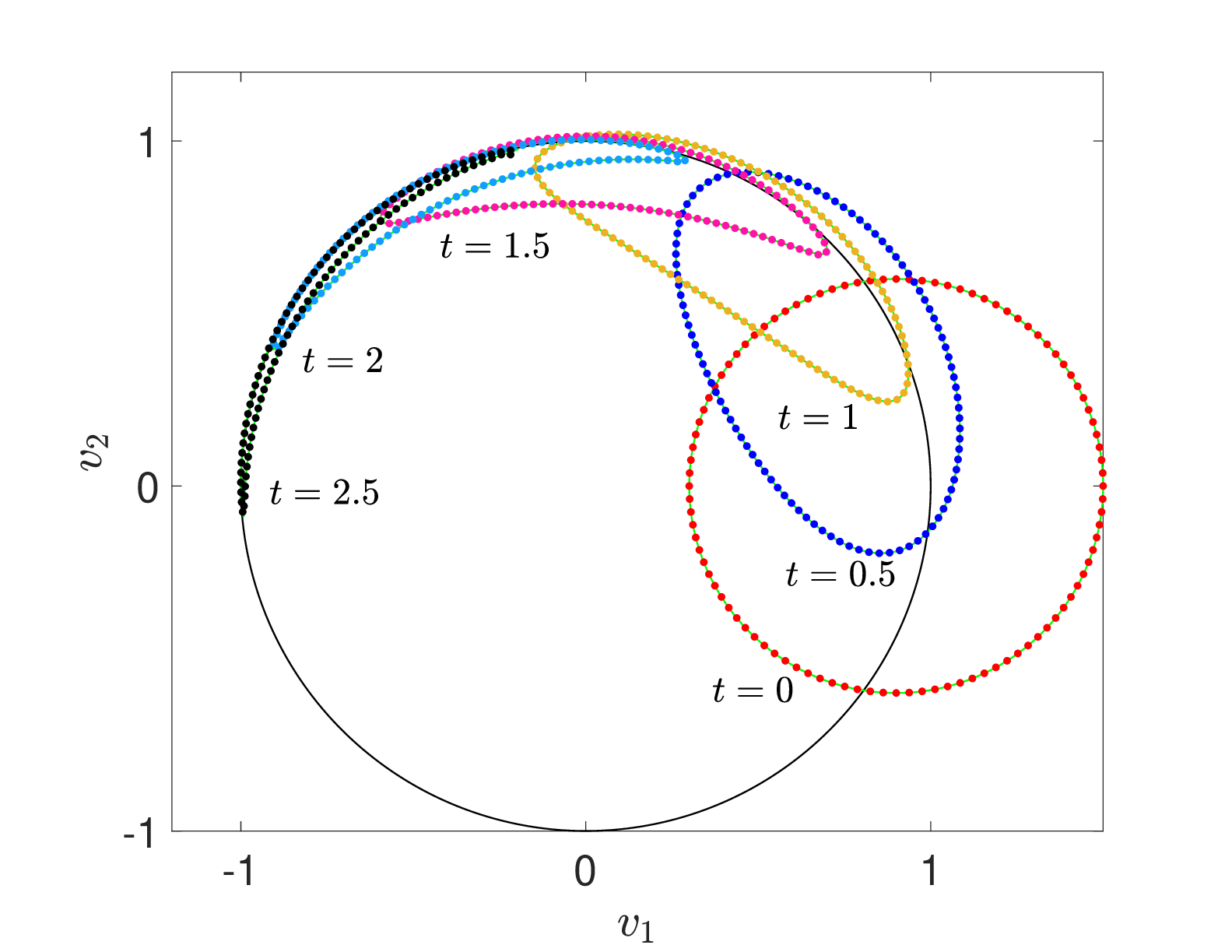}\\
     (a)    & (b)
    \end{tabular}
    
    \caption{Evolution of a closed curve of initial values for the two-dimensional
  dynamical system \eqref{normHopf};  plot of $6$ curves at times
  $t=0,0.5,1, 1.5, 2$ and $2.5$ in a space-time diagram (left)
  and overlapping in a space diagram (right); see Figure \ref{fig5.6}
  for more details on the numerical computation of the dotted
lines in the right picture.}
    \label{fig1}
\end{figure}

To illustrate the procedure we consider 
the simplest case of all, namely closed curves which evolve under
the flow of a dynamical system and for which a parameterization is determined
during the computation. Concretely, the state space is 
$ C^{\infty}(S^1,\R^d)$  where $ S^1$ is the unit circle in $\R^2$ and
$F$ is a Nemitzky operator, generated by a vector field $f:\R^d \to \R^d$,
i.e.\
\begin{equation} \label{e1d}
  F(u)(x)= f(u(x)), \quad x\in S^1, u\in C^{\infty}(S^1,\R^d).
\end{equation}
Figure \ref{fig1} shows the evolution of an initial circle for the two-dimensional
dynamical system $\dot{u}=f(u)$ representing the normal form of the Hopf bifurcation \eqref{normHopf}.
The initial circle of radius $0.6$ with center at $(0.9,0)$
intersects the periodic orbit of the system  which happens to be the unit circle. 
Then,
the dynamics of the system
stretches the initial curve into a long, banana-like closed curve.  This dynamics is
followed by discretizing \eqref{e1b},\eqref{e1c} with a suitable phase
condition, see Section \ref{sec4.2} for more details.

Of course, one can obtain a similar picture by solving the initial value problem
$\dot{u}=f(u)$ for a number of points on the initial circle and then
connecting them by interpolation. However, as we will see, if the system's dynamics
drives the points apart/closer, the shape
of the curve is quickly lost and interpolation becomes useless. Therefore, 
one aim  of our reformulation is to preserve the connectedness of
the initial manifold and to adapt the coordinate systems {\em during integration}
so that the manifold's evolution becomes visually clear.
We achieve this  by incorporating not only the time derivative, but also spatial derivatives of the solution into an {\it energy-based functional}, which is then  minimized 
during time integration.

In case of the Nemitzky operator \eqref{e1d} the additional term
$\partial_x v \mu$ turns the Freezing system 
\eqref{e1b},\eqref{e1c} into a first order hyperbolic system which we solve
by the Lax-Wendroff method. If the operator $F$ involves higher order spatial
derivatives, for example the Laplace-Beltrami operator, then the additional term is of lower order and the type of the system remains unchanged.
Applications of this type will be briefly discussed in the Conclusion section~\ref{sec6}.


\sect{Diffeomorphism equivariance and the Freezing system}
\label{sec1}
In this section we define the abstract setting of spaces and operators. Our main aim is to derive the Freezing
system \eqref{e1c}, \eqref{e1d} in a formal manner. 
By ``formal" we mean that we
do not yet specify the topology on the respective spaces nor  the smoothness properties
of solutions and operators.
  We will also not detail the analysis on the infinite dimensional Lie group
  $\mathrm{Diff}(M)$ and of its action on the phase space.
  In \ref{sec2.3} we will indicate how all this can be
  done within the rigorous framework of global analysis with convenient
  vector spaces developed in \cite{KM97}. A more extensive account of
  embedding our results into this abstract concept is
  laid out in \cite{Roe19}.

\subsection{The setting and notations}
\label{sec1.1}
In the following let $M$ be a real compact $p$-dimensional orientable
$C^{\infty}$-manifold with tangent bundle $TM$. 
For some $d\in \N$ our phase space is 
\begin{equation} \label{e2}
\mathcal{U}_d= C^{\infty}(M,\R^d),
\end{equation}
with the space of $C^{\infty}$-vector fields on $M$ denoted by
\begin{equation} \label{e2a}
  \mathcal{V}= C^{\infty}(M,TM),
\end{equation}
and the infinite-dimensional  group of diffeomorphisms of M denoted by
\begin{equation} \label{e2b}
  G= \mathrm{Diff}(M).
 \end{equation}
  

As noted in the Introduction, whenever we speak of spatial derivatives we mean derivatives with respect to the coordinates on $M$.
\begin{remark} \label{rem1}
We use the normed linear space $\R^d$ as an image space which turns
$\mathcal{U}_d$ into a vector space. It may be given the structure of a
locally convex topological
vector space; see \ref{sec2.3}. This will simplify our
considerations. It is possible
to replace $\R^d$ by a $d$-dimensional
$C^{\infty}$-manifold with a Riemannian structure, but we avoid the ensuing
complications.
\end{remark}
The right-hand side of \eqref{e1} is defined by a vector field which maps $\mathcal{U}_d$ into its tangent bundle $T \mathcal{U}_d$, i.e.
\begin{equation*} \label{e9}
F:\mathcal{U}_d \rightarrow T \mathcal{U}_d, \quad \textrm{ where } F(u) \in T_u \mathcal{U}_d \quad
\text{for} \; u \in \mathcal{U}_d.
\end{equation*}
Since $\mathcal{U}_d$ is a linear space we can identify for each $u \in \mathcal{U}_d$ its tangent
space $T_u \mathcal{U}_d=\{u\}\times \mathcal{U}_d$ with $\mathcal{U}_d$ itself.
The image $F(u)$ may involve spatial derivatives
of $u$ with respect to $x$ in $M$ of any order, but the  examples mentioned in the introduction involve
only derivatives of order at most two. For our model example from the introduction  we have
 $M=S^1$ and $F$ is given by a Nemitzky operator \eqref{e1d},\ i.e. it
involves no spatial derivatives.

The underlying time interval is $\mathcal{I}=[0,T)$ where either $T>0$
  or $T=\infty$. The image of any  function
\begin{align} \label{e3a}
  u: \begin{cases} \mathcal{I} & \to \quad \mathcal{U}_d, \\
      t & \mapsto \quad u(t)
  \end{cases}
\end{align}
is a function of $x \in M$ itself. In accordance
with standard PDE notation we may view $u$ either as a mapping into a
function space as  in \eqref{e3a} or as a function
$u:\mathcal{I}\times M \to \R^d$  and write
\begin{align} \label{e3b}
  [u(t)](x)= u(t,x) \in \R^d, \quad t \in \mathcal{I}, x \in M.
\end{align}
By suppressing the $x$-dependence we abbreviate  this as
$u(t) =u(t,\cdot)$, $ t \in \mathcal{I}$. When including smoothness properties
  of $u$, the identification \eqref{e3b}  is already a subtle issue,
  called the exponential law in global analysis; see \cite[Ch.I,IX]{KM97}
  and the discussion in \ref{sec2.3}.
  
Correspondingly, the derivative $\frac{d}{dt}u(t) \in T_{u(t)}\mathcal{U}_d=\mathcal{U}_d$ of the function \eqref{e3a} may be written in two equivalent forms 
\begin{align} \label{e3c}
  \left[ \frac{d}{dt}u(t)\right](x) = \partial_tu(t,x) , \quad t \in \mathcal{I}, x \in M,
\end{align}
where $\partial_tu(t,x) \in \mathbb{R}^d $ denotes the derivative (in time) of
the curve $t \mapsto u(t,x) \in  \mathbb{R}^d $ for fixed $x \in M$.
With these identifications the initial value problem \eqref{e1} reads pointwise
\begin{align*}
  \partial_tu(t,x) = F(u(t,\cdot))(x) , \quad  t \in \mathcal{I}, x \in M.
\end{align*}

Our aim is to consider equivariance with respect to
the group $G$  of diffeomorphisms \eqref{e2b}.
For $g\in G$ we denote its tangent map  (or total derivative) by 
\begin{equation} \label{e5}
T_xg:TM \rightarrow TM.
\end{equation}
It maps the tangent bundle into itself such that for each
$x \in M$ the map
\begin{equation*} \label{e5a}
  T_xg(x):T_x M \rightarrow T_{g(x)}M 
\end{equation*}
is a linear homeomorphism.

The group operation in $G$ is composition $g_1\circ g_2$, $g_1,g_2\in G$
with the unit element given by the identity $\one$, $\one(x) = x, x \in M$,
and the inverse of $g \in G$ denoted by $g^{-1}$.
With this, $G$ becomes an infinite dimensional Lie group. The analysis
on such manifolds is rather involved; see the monograph \cite{KM97}
and the discussion in \ref{sec2.3}.

The group $G$ acts on $\mathcal{U}_d$  by composition, that is
\begin{equation}\label{e6}
a:G\times \mathcal{U}_d \rightarrow \mathcal{U}_d, \quad
a(g,u)= u \circ g^{-1}, \quad g\in G, u\in \mathcal{U}_d.
\end{equation}
Note that $a(g,u)$ is linear with respect to $u \in \mathcal{U}_d$ but
nonlinear with respect to $g \in G$, in general.

If we have a time-dependent family of diffeomorphisms (cf. \eqref{e3a})
\begin{align} \label{e6a}
  g: \begin{cases} \mathcal{I} & \to \quad G, \\
      t & \mapsto \quad g(t),
  \end{cases}
\end{align}
then we use the same identification as in \eqref{e3b}, i.e.
\begin{align} \label{e6b}
  [g(t)](x)= g(t,x) \in M, \quad t \in \mathcal{I}, x \in M.
\end{align}
Taking the (formal) $t$-derivative of $g$ in \eqref{e6a}, \eqref{e6b} then leads to the equation
(cf. \eqref{e3c})
\begin{align} \label{e6c}
    \left[ \frac{d}{dt} g(t) \right](x) = \partial_t g(t,x) \in T_{g(t,x)}M,
    \quad x \in M.
\end{align}
Thus  we find that for every $t \in \mathcal{I}$ the map $ \frac{d}{dt}g(t):M \to TM$ is an element of the space $\mathcal{V}$ of vector fields, defined by the relation \eqref{e6c}.
  Since $g$ depends on $t$ we now write the $x$-derivative \eqref{e5} as
  a partial derivative,
  i.e.
  \begin{align} \label{e6d}
    \partial_xg(t,\cdot):TM \to TM, \quad \partial_xg(t,x):T_xM \to T_{g(t,x)}M.
  \end{align}
  
Our main assumption is that the vector field $F$ is equivariant
with respect to the action $a$ of the group $G$, i.e.
\begin{equation} \label{e10}
F(a(g,u))= a(g,F(u)), \quad g\in G, u \in \mathcal{U}_d.
\end{equation}
Note that both sides are elements of $\mathcal{U}_d$ since $\mathcal{U}_d$ is
linear space. Equation \eqref{e10} reads explicitly
$F(u\circ g^{-1}) = (F(u))\circ g^{-1}$ for all $g \in  G$,  or by setting
$\gamma=g^{-1}$,
\begin{align} \label{e11}
  F(u\circ \gamma) = (F(u))\circ \gamma \quad \text{for all}
  \quad \gamma \in G, u \in \mathcal{U}_d.
\end{align}

\begin{remark} \label{rem2}
  If the image space $\R^d$ is a manifold as in Remark \ref{rem1},  then \eqref{e10} should be replaced
  by $F(a(g,u))=T_ua(g,F(u))$ where $T_ua$ denotes the tangent map
  of the action with respect to the second variable.
 \end{remark}

\begin{example} \label{ex2}
  As mentioned in the introduction,
  the simplest case of equivariance occurs if
  $F(u)$ does not depend on spatial derivatives of $u$, i.e. if $F$ is a Nemitzky operator \eqref{e1d}.
  Indeed, we then have for $g\in G$, $u\in \mathcal{U}_d$
\begin{equation} \label{e12}
F(u\circ g^{-1})(x)= f(u(g^{-1}(x))) = [F(u) \circ g^{-1}](x), \quad x \in M.
\end{equation}
 We imagine that a vector field
$f: \R^d \rightarrow \R^d$ is given, and define
$F$  by \eqref{e1d}. Then the differential equation
\eqref{e1} describes the motion of the initial manifold $u_0(M)$ induced by the 
flow of the dynamical system
\begin{equation} \label{eq:dynsyst}
  \dot{v}=f(v).
\end{equation}
Our aim is to endow $M$ with a suitable coordinate system 
with favorable properties \underline{still to be determined}. The particular case 
$M=S^1$ of flowing closed curves will be 
the test case for the numerical methods.
\end{example}

\subsection{Derivation of the Freezing system}
\label{sec2.1}
Following the approaches in \cite{BT04}, \cite{RKML03} we extend the method of
Freezing to evolution equations that are equivariant with respect
to the infinite dimensional Lie group of diffeomorphisms
discussed above.

The idea is to write the solution of \eqref{e1} as the group action
\eqref{e6} of an unknown
element $g(t)\in G,t\ge 0$ applied to a just as unknown function $v(t) \in
\mathcal{U}_d$, i.e.
\begin{equation*} \label{e34}
u(t)=a(g(t),v(t)), \quad t \in \mathcal{I}.
\end{equation*}
Inverting $g(t,\cdot)$ and using the notation $u(t,\cdot)$ above leads to
\begin{equation} \label{e12a}
  v(t,\cdot) = u(t,g(t,\cdot)), \quad t\in \mathcal{I}.
\end{equation}
Differentiating \eqref{e12a} with respect to $x \in M$ gives
\begin{equation} \label{e12b}
  \partial_x v(t,\cdot) = \partial_x u(t,g(t,\cdot)) \partial_xg(t,\cdot).
\end{equation}
Note that $\partial_xg(t,\cdot)$ has a smooth inverse as in \eqref{e6c}
since $g$ is a smooth diffeomorphim.
Differentiating \eqref{e12a} with respect to $t \in \mathcal{I}$ and using
equation \eqref{e1} as well as
equivariance \eqref{e11} and \eqref{e12a}, \eqref{e12b} yields
\begin{align*}
  \partial_tv(t,\cdot) = &  \partial_x u (t,g(t,\cdot))\partial_tg(t,\cdot)
  + \partial_t u (t, g(t,\cdot)) \\
  =& \partial_x u (t,g(t,\cdot))\partial_tg(t,\cdot) +
  \left[F(u(t,\cdot))\right](g(t,\cdot))\\
  = &\partial_xv(t,\cdot) (\partial_xg(t,\cdot))^{-1} \partial_t g(t,\cdot)
  +F(u(t,g(t,\cdot))) \\
  = &\partial_xv(t,\cdot) (\partial_xg(t,\cdot))^{-1} \partial_t g(t,\cdot)
  +F(v(t,\cdot)).
\end{align*}
This equation suggests to introduce the vector field
\begin{equation*} \label{e12c}
  \mu(t,\cdot)= (\partial_xg(t,\cdot))^{-1} \partial_tg(t,\cdot): M
  \to TM,
\end{equation*}
where $\mu(t,x)=(\partial_xg(t,x))^{-1} \partial_tg(t,x) \in T_xM$ for $x \in M$ follows from \eqref{e6c} and \eqref{e6d}.

Then the Cauchy problem \eqref{e1} is equivalent to finding solutions
$v :\mathcal{I} \to \mathcal{U}_d$, $g: \mathcal{I} \to G$ and a time-dependent
    vector field  $\mu:\mathcal{I} \to \mathcal{V}$ which solve the
      {\it Freezing system}
\begin{align} \label{e37a}
  \partial_tv(t,\cdot)= & F(v(t,\cdot))+ \partial_xv(t,\cdot) \mu(t,\cdot) ,\quad v(0,\cdot)=u_0(\cdot), \\
  \label{e37b}
\partial_t g(t,\cdot) = & \partial_xg(t,\cdot) \mu(t,\cdot), \quad g(0,\cdot)= \one. 
\end{align}
Writing this pointwise for $x\in M$ we obtain a first order system for the functions
\begin{equation} \label{e37c}
\begin{aligned}
 v:&\; \mathcal{I}\times M \rightarrow \R^d, \quad (t,x) \mapsto v(t,x),\\
 g:&\; \mathcal{I}\times M \rightarrow M, \quad (t,x) \mapsto g(t,x),\\
\mu:& \mathcal{I} \times M \rightarrow TM, \quad (t,x) \mapsto
\mu(t,x) \in T_x M.
\end{aligned}
\end{equation}
The type of the PDE \eqref{e37a} depends on the number of derivatives involved in the vector field $F$. If $F$ derives from a Nemitzky operator as in \eqref{e1d}, then the first order term $\partial_x v \mu$ in
\eqref{e37a} is dominant, and the system is  of symmetric hyperbolic type. If $F$ involves second
order derivatives then the extra terms $\partial_x v \mu$ are generally of lower order
and will not change the type of $F$.
Following \cite{RKML03}, equation \eqref{e37b} is called the {\it reconstruction equation},
since it allows to determine the diffeomorphisms $g(t,\cdot)$ {\em a-posteriori},
once the solutions $v$ and $\mu$ have been computed.
Because of the extra unknown vector field $\mu$ the system is not yet closed.
In  Section \ref{sec2} we discuss several possibilities for
closing the system by suitable phase conditions.

\sect{Phase conditions}
\label{sec2}

\subsection{Minimization principles} \label{sec2.2}
Phase conditions are usually derived from optimization principles
which aim at efficient numerical computations and at understanding the longtime dynamics of a system. From now on we assume  $M$ to be a Riemannian
manifold and a submanifold of some $\R^n$. Further, we make use of the inner product
 on $\mathcal{U}_d$  given by
\begin{equation*} \label{e30}
( u,v )_{L^2(M,\R^d)}= \int_{M} ( u(x),v(x) )_{\R^d} \
\mathrm{do}(x), \quad u,v \in \mathcal{U}_d,
\end{equation*}
where $( \cdot,\cdot )_{\R^d}$ denotes the Euclidean inner product
in $\R^d$ and the integral is defined with respect to the volume element induced by the metric tensor.

The phase conditions which close the system \eqref{e37a}, \eqref{e37b}
are of the general form
\begin{equation} \label{e31}
  \psi(v,\mu,g) = 0,
\end{equation}
with a given {\it phase function}
\begin{equation} \label{e32}
  \psi:\mathcal{U}_d \times \mathcal{V} \times G \to C^{\infty}(TM,\R)=:\cZ.
\end{equation}
Altogether, the system \eqref{e37a},\eqref{e37b} with the side constraint
\eqref{e31} constitutes a partial differential algebraic equation
(PDAE) which is to be solved for the functions $v,\mu,g$ from \eqref{e37c}.
If the phase function $\psi$
from \eqref{e32} does not involve $g$, then one may first solve the
reduced PDAE \eqref{e37a}, \eqref{e31} and afterwards solve the reconstruction
equation \eqref{e37b}.

Each of the following phase conditions is obtained from the necessary first
order conditions for a local minimum of a {\it phase functional}
\begin{equation*} \label{e33}
  \rho:\mathcal{U}_d \times \mathcal{V} \times G \to \R, \quad
  (w,\nu,\gamma) \mapsto \rho(w,\nu,\gamma).
\end{equation*}
We distinguish two types of phase conditions based on minimizing either
with respect to $\nu \in \mathcal{V}$ for $v,g$ fixed, or with respect to $\gamma \in G$ for $v,\mu$ fixed :
\begin{equation*}
\begin{aligned}
\mu& = \mathrm{argmin}_{\nu \in \mathcal{V}}\rho(v,\nu,g) \quad \text{type I}, \\
  \one & = \mathrm{argmin}_{\gamma \in G}\rho(v,\mu,\gamma) \quad \text{type II}.
  \end{aligned}
  \end{equation*}
We consider four phase functionals
defined for $(w,\nu,\gamma)\in \mathcal{U}_d \times \mathcal{V} \times G$
as follows:
\begin{enumerate}
\item[(1)] \textbf{orthogonal phase functional (type I)}
  \begin{equation*} \label{e35b}
    \rho_1(w,\nu,\gamma) = \|F(w)+ \partial_x w \nu \|^2_{L^2(M,\R^d)}.
  \end{equation*}
  \item[(2)] \textbf{fixed phase functional (type II)}
  \begin{equation*} \label{e35a}
    \rho_2(w,\nu,\gamma)= \|w - \hat{v}\circ \gamma \|^2_{L^2(M,\R^d)},
  \end{equation*}
  where $\hat{v}\in \mathcal{U}_d$ is a given template function.
  \item[(3)] \textbf{energy based phase functional (type II)}
    \begin{equation*} \label{e35c}
      \rho_3(w,\nu,\gamma)=\|\partial_x(w \circ \gamma)\|^2_{L^2(M;TM, \R^d)},
    \end{equation*}
    which is the integral of the Frobenius (Hilbert-Schmidt) norm
    of the linear map $\partial_x w(x):T_xM \to \R^d$ with $T_xM$ considered
    as a subspace of $\R^n$:
    \begin{equation*}
      \|\partial_x w\|^2_{L^2(M;TM,\R^d)} = \int_M \| \partial_x w(x)\|_F^2
      \mathrm{do}(x)= \int_M \sum_{\ell=1}^d \|\mathrm{grad}w_{\ell}(x)\|^2_{\R^n} \mathrm{do}(x).
    \end{equation*}
    We refer to \cite[Ch.XI]{AEIII09} for a representation of the gradient
    in local coordinates.
  \item[(4)] \textbf{semi-discrete energy based phase functional (type I)}
    \begin{equation*}\label{e35d}
      \rho_4(w,\nu,\gamma)= \| \partial_x w + \Delta t \partial_x(F(w)+\partial_x w \nu)\|^2_{L^2(M;TM,\R^d)},
    \end{equation*}
    where $\Delta t$ is a small parameter, taken to be the time step of the
    numerical scheme.
\end{enumerate}
Our motivation behind these choices is as follows:
The first two functionals are familiar from  the finite
dimensional case; see \cite{BOR14,BT04,BT07}.  The functional $\rho_1$ 
aims at minimizing the norm of
$\partial_t v=F(v)+\partial_x v \mu$ (see \eqref{e37a}),
i.e.\ one tries to let $v$ move as little as possible.
The term ``orthogonal phase functional" is suggested by the resulting
phase condition \eqref{eq:orthphase2} derived below.
The functional $\rho_2$ uses an (arbitrary) template function $\hat{v}$
which may for example be taken as the initial data $\hat{v}=u_0$ or as $\hat{v}=v(\hat{t})$ at some later time $\hat{t}>0$.
Minimizing $\rho_2$ w.r.t. $\gamma$ then requires $\hat{v}$ to be the closest point to $v(t)$ on the group orbit of $\hat{v}$. Thus one tries
to keep the solution as close as possible to $\hat{v}$ during
integration.

The two functionals $\rho_3,\rho_4$ form a new selection since they involve the spatial derivative
of the solution $v$. While the first two phase conditions aim at simplifying
the longtime dynamics of the solution, the new conditions aim at avoiding
large gradients of the solution during integration. Such a goal is suggested
by the numerical experiments shown in Section \ref{sec:results}. It turns out
that $\rho_3$ can lead to inconsistent initial data for the resulting
PDAE and thus may cause numerical difficulties; see Section \ref{sec4.3}.
The modified functional $\rho_4$ alleviates this problem
by realizing a steepest descent of the energy norm within some
time-step $\Delta t>0$. It is based on the truncated  Taylor expansion
\begin{align*}
  \partial_xv(t+\Delta t,\cdot)& \approx \partial_xv(t,\cdot)
  + \Delta t \partial_t \partial_xv(t,\cdot)= \partial_xv(t,\cdot)
  + \Delta t \partial_x \partial_tv(t,\cdot)\\
  & = \partial_xv(t,\cdot) + \Delta t \partial_x\big(F(v(t,\cdot))+
  \partial_xv(t,\cdot)\mu(t,\cdot)\big).
\end{align*}
In principle, we can use any small parameter instead of $\Delta t$, but
using the time-step $\Delta t$ of the numerical scheme is a very natural
choice and works well in practice; see Section \ref{sec4}.

\subsection{Derivation of phase conditions}
\label{sec3.2}
In this subsection we set up the necessary conditions for local minima
of the phase functionals from Section \ref{sec2.2}. As we will see,
the type numbers I and II correspond to the (differentiation) index of
the resulting PDAE (for a survey of index concepts for DAEs see e.g. 
\cite[Ch.1.2]{KM06}).

As usual, the necessary conditions are derived from the  fundamental lemma
of the calculus of variations. We use the following version  (see
\cite[A.2]{LLM18} and \cite[Section 3.2]{Roe19} for a proof of this specific case):
\begin{lemma}\label{lem3.1}
  Let $\varphi\in \cZ$ be a smooth functional on the tangent bundle
  satisfying
  \begin{equation} \label{eq3:fundlem}
    \begin{aligned}
      \varphi(x,\lambda X) & = \lambda \varphi(x,X), \quad \forall x\in M,
      X\in T_xM, \lambda \in \R,\\
      \int_{M}\varphi(\nu(x)) \mathrm{do}(x)& = 0,
      \quad \forall \nu \in \mathcal{V}.
    \end{aligned}
  \end{equation}
  Then $\varphi(x,X)=0$ holds for all $x \in M$, $X \in T_xM$.
\end{lemma}
\begin{remark} Formally, the tangent space $T_xM$ at $x \in M$
is defined as the set of pairs $(x,X)$ where $X$ lies in a
$p$-dimensional linear space. For simplicity, we often write
 $X \in T_xM$ instead of  $(x,X)\in T_xM$; see \eqref{eq3:fundlem}.
\end{remark}

Let us apply Lemma \ref{lem3.1} to the four phase functionals above. 
\begin{enumerate}
\item[(1)] {\it Orthogonal phase condition (type I)}\\
  If $\rho_1$  has a local minimum at $\mu$
  then  for each $\nu \in \mathcal{V}$ it holds that
\begin{equation*}
  \label{eortho}
  \begin{aligned}
    0 &= \frac{\mathrm{d}}{\mathrm{d \tau}}\rho_1(w, \mu+\tau \nu,\gamma)_{|\tau=0}
    =\Big[ \frac{\mathrm{d}}{\mathrm{d \tau}}\|F(w)+\partial_x w (\mu + \tau
      \nu)\|^2_{L^2(M, \R^d)}\Big]_{|\tau=0}\\
      &= 2 \int_M(F(w)(x)+ \partial_xw(x)\mu(x), \partial_x w(x) \nu(x))_{\R^d}
      \mathrm{do}(x).
\end{aligned}
\end{equation*}
Thus Lemma \ref{lem3.1} applies to
\begin{align*}
  \varphi(x,X)=(F(w)(x)+ \partial_xw(x)\mu(x), \partial_x w(x) X)_{\R^d},
  \quad x\in M, X \in T_xM
\end{align*}
and yields  a phase condition of the form \eqref{e31}
\begin{equation} \label{eq:orthphase}
  \begin{aligned}
    \psi_1(v,\mu,g)(x,X)&=(F(v)(x)+ \partial_xv(x)\mu(x), \partial_x v(x) X)_{\R^d}=0,\\
    & \forall x\in M, X\in T_xM.
    \end{aligned}
\end{equation}
In the time dependent setting of \eqref{e37a} this may be written as
\begin{equation} \label{eq:orthphase2}
  \begin{aligned}
   (\partial_t v(t,x), \partial_x v(t,x) X)_{\R^d}=0, \quad
    & \forall x\in M, X\in T_xM,
    \end{aligned}
\end{equation}
i.e.  {\em time and space derivatives of $v$ are orthogonal}.
Figure \ref{fig2} gives a geometric interpretation of the phase condition
\eqref{eq:orthphase} for Example \ref{ex2}  with a two-dimensional
manifold $M$. The term $v_x \mu$ is added to the vector field $f$ in
\eqref{e37a} in such a way that the resulting vector field is orthogonal to
the current manifold $v(t,M)$ of $v$-images.

\begin{figure}[!ht] \label{fig2}
  \includegraphics[width=\textwidth]{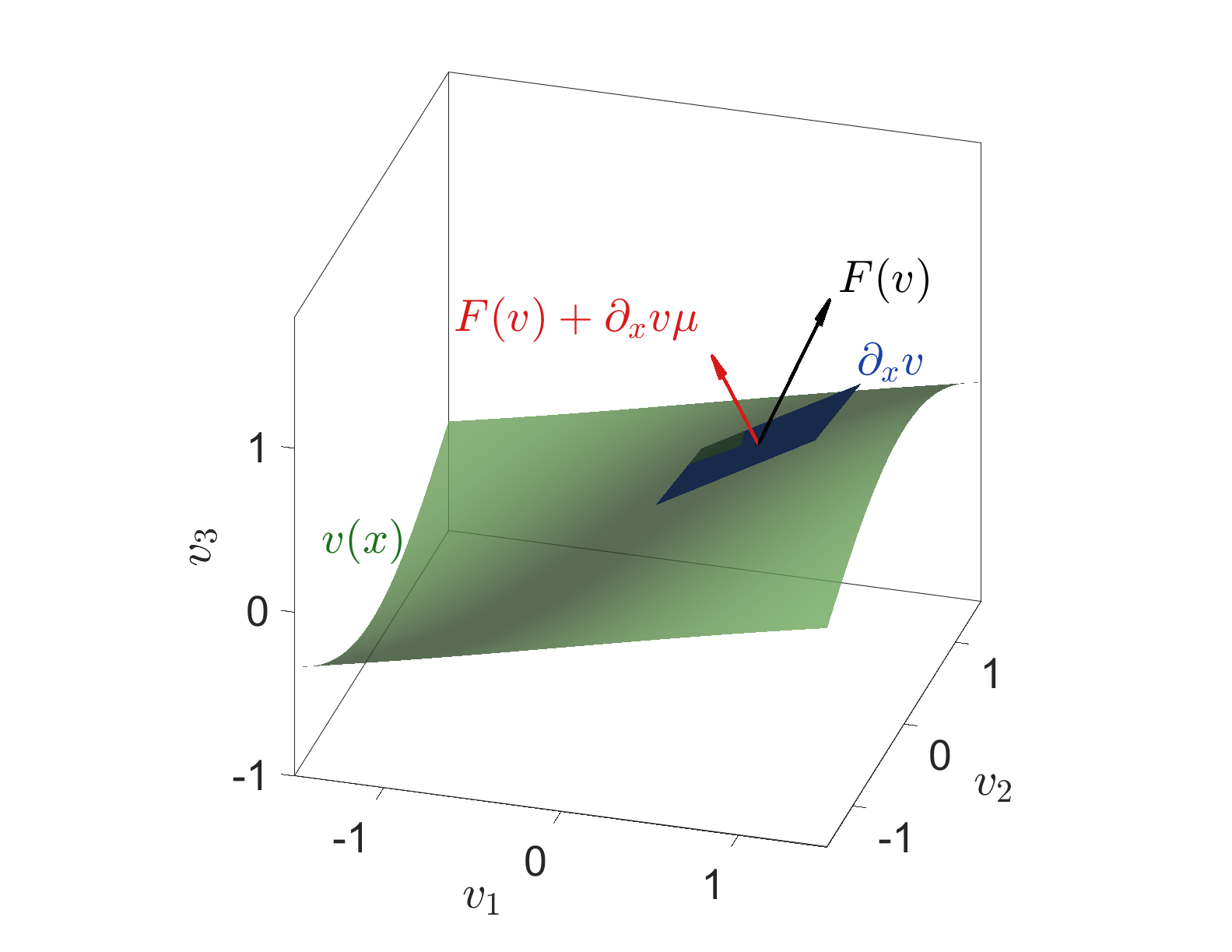}
\caption{Illustration of the orthogonal phase condition \eqref{eq:orthphase}}
\end{figure}

Selecting a basis in $T_xM$ one can represent $\partial_xv(t,x)\in \R^{d,p}$ as
a $d \times p$ matrix and rewrite \eqref{eq:orthphase} as the linear system
\begin{equation} \label{eq:mateq}
  \partial_xv(t,x)^{\top}\partial_xv(t,x) \mu(t,x)= - \partial_xv(t,x)^{\top}F(v(t,\cdot))(x).
\end{equation}
If $\partial_xv(t,x)$ has maximum rank  $p$ for all $x \in M$ then this system has a unique
solution $\mu(t,x)$, which means that the overall PDAE has (differentiation) index $1$. One can also use \eqref{eq:mateq} at $t=0$ with $v(0,x)=u_0(x)$ to determine consistent initial values for $\mu(0,x)$, $x\in M$.

\item[(2)] {\it Fixed phase condition (type II)} \\
  The variation of the diffeomorphism proceeds as follows. For any
  vector field $\nu \in \mathcal{V}$ there exists a smooth path $\gamma(\tau,\cdot)\in G$, $|\tau|< \varepsilon$ for some $\varepsilon >0$ such that $\gamma(0,\cdot)=\one$ and
  $\nu(x)=\partial_{\tau}\gamma(0,x)\in T_xM$ for all $x\in M$; see the
  remarks on the tangent bundle of $G$ in Appendix \ref{sec2.3}.
  Then we obtain
  \begin{equation*} \label{efixed}
    \begin{aligned}
    0& = \frac{\mathrm{d}}{\mathrm{d \tau}} \rho_2(w,\nu,\gamma(\tau,\cdot))_{|\tau=0}
    =\frac{\mathrm{d}}{\mathrm{d \tau}}\left[ \|w-\hat{v}(\gamma(\tau,\cdot)))\|^2_{L^2(M,\R^d)}\right]_{|\tau=0}\\
    & = - 2 \int_M\big(w(x)-\hat{v}(\gamma(\tau,x)),
    \partial_x\hat{v}(\gamma(\tau,x))\partial_{\tau}\gamma(\tau,x)\big)_{\R^d}
    \mathrm{do}(x)_{|\tau=0}\\
    &= - 2 \int_M \big(w(x) - \hat{v}(x),\partial_x \hat{v}(x) \nu(x)\big)_{\R^d} \mathrm{do}(x).
  \end{aligned}
  \end{equation*}
  Thus Lemma \ref{lem3.1} applies to
  \begin{align*}
    \varphi(x,X)= \big( w(x) -\hat{v}(x), \partial_x \hat{v}(x)X\big)_{\R^d},
    \quad x \in M, X \in T_xM
    \end{align*}
  and yields for the solution $v(t,\cdot)$ at time $t$ the algebraic condition
  \begin{equation} \label{eq:fixedphase}
    0= \big( v(t,x) -\hat{v}(x), \partial_x \hat{v}(x)X\big)_{\R^d}
    \quad \forall x \in M, X \in T_xM.
  \end{equation}
  This condition does not yet involve the algebraic variable $\mu$ so that
  the PDAE is of higher index. As usual for DAEs, we differentiate the
  constraint \eqref{eq:fixedphase} w.r.t. $t$ and use \eqref{e37a} to
  obtain for $x \in M$, $X\in T_xM$ the hidden constraint
  \begin{equation} \label{eq:hidden}
    \begin{aligned}
      0&= \big( \partial_tv(t,x),\partial_x\hat{v}(x) X \big)_{\R^d} \\
    &=\big(F(v(t,\cdot))(x)+ \partial_xv(t,x)\mu(t,x), \partial_x\hat{v}(x)X
    \big)_{\R^d} \\
    & =: \psi_2(v(t,\cdot),\mu(t,\cdot),g)(x,X).
    \end{aligned}
  \end{equation}
  As above this equation is uniquely solvable for $\mu(t,x)$
  if $\partial_xv(t,x)^{\top}\partial_x \hat{v}(x)$ is invertible, in which
  case we arrive at a PDAE of index $2$. Note that \eqref{eq:fixedphase} is automatically satisfied at $t=0$ for the choice $\hat{v}=u_0$. Together with
  \eqref{eq:hidden} this implies that \eqref{eq:fixedphase} holds for
  all $t \in \mathcal{I}$.
\item[(3)] {\it Energy based phase condition (type II)}\\
  Proceeding as in the previous item we find the necessary
  condition
  \begin{equation*} \label{eq:energybased}
    \begin{aligned}
      0 & = \frac{\mathrm{d}}{\mathrm{d \tau}}\rho_3(w,\nu, \gamma(\tau,\cdot))_{|\tau=0}\\
      & = 2 \sum_{\ell=1}^d \int_M \big(\mathrm{grad}[\partial_{\tau}(w_{\ell}\circ \gamma(\tau,\cdot))],\mathrm{grad}[w_{\ell}\circ \gamma(\tau,\cdot)]\big)_{\R^n} \mathrm{do}_{|\tau=0}\\
      & = 2 \sum_{\ell=1}^d\int_M \Big( \mathrm{grad}[\partial_xw_{\ell} \circ \nu](x),\mathrm{grad}[w_{\ell}](x) \big)_{\R^n} \mathrm{do}(x) \\
      & =- 2 \sum_{\ell=1}^d \int_M [\partial_xw_{\ell}(x) \nu(x)] \Delta_Mw_{\ell}(x) \mathrm{do}(x).
    \end{aligned}
  \end{equation*}
  In the last step we used integration by parts and the Laplace-Beltrami operator
  $\Delta_M$ on $M$; see \cite[Ch.XI]{AEIII09}
  and Section \ref{sec6} for a representation of $\Delta_M$
  in local coordinates. Now Lemma \ref{lem3.1} applies to
  \begin{align*}
    \varphi(x,X)= \sum_{\ell=1}^d [\partial_xw_{\ell}(x) X] \Delta_Mw_{\ell}(x)
    = \big(\partial_xw(x)X, \Delta_Mw(x)\big)_{\R^d}
  \end{align*}
  where $\partial_xw X$ and $\Delta_M w$ are defined componentwise.
  We obtain the phase condition
  \begin{equation} \label{eq:phaseenergy}
    0 = \big( \partial_xv(t,x)X, \Delta_Mv(t,x)\big)_{\R^d} \quad \forall x\in M, X \in T_xM,
  \end{equation}
  which we differentiate w.r.t. time to find the hidden constraint
  \begin{equation} \label{eq:phaseenergy2}
    \begin{aligned}
      0& = \psi_3(v(t,\cdot),\mu(t,\cdot),g)(x,X) \\
      & = \Big( \partial_x\big(F(v(t,\cdot))+ \partial_xv(t,\cdot)\mu(t,\cdot)\big)(x)X,
      \Delta_Mv(t,x) \Big)_{\R^d} \\
      & + \Big( \partial_xv(t,\cdot)X, \Delta_M[F(v(t,\cdot))+ \partial_x v(t,\cdot) \mu(t,\cdot)](x) \Big)_{\R^d} \; \forall x\in M, X \in T_xM.
    \end{aligned}
  \end{equation}
  Integration by parts for the first term and evaluation of $\Delta_M$ in the second term
  shows that \eqref{eq:phaseenergy2} is a linear second order
  elliptic problem for $\mu(t,\cdot)$ of the form:
  \begin{align} \label{eq:elliptic}
    \big[\partial_xv(t,\cdot)^{\top}\partial_xv(t,\cdot)\big] \Delta_M \mu(t,\cdot)
    = a(t,\cdot)\mu(t,\cdot)+ b(t,\cdot),
  \end{align}
  where the coefficients $a$ and $b$ depend on spatial derivatives of $v(t,\cdot)$
  up to order $3$.
Moreover, let us note that we cannot assume 
 the condition \eqref{eq:phaseenergy} to hold for the initial
 data $v(0,\cdot)=u_0$, in general. Therefore, contrary to the fixed
 phase condition, we also cannot conclude from \eqref{eq:elliptic} that
 \eqref{eq:phaseenergy} holds for all $t \in \mathcal{I}$.
\item[(4)] {\it semi-discrete energy based functions (type I)} \\
  The analysis proceeds along the lines of the orthogonal phase condition
  and leads to
  \begin{equation} \label{eq:semidiscrete}
    \begin{aligned}
      0& = \psi_4(v(t,\cdot),\mu(t,\cdot),g)(x,X) \\
      & = \Big( \partial_xv(t,x) X, \Delta_M\big(v(t,\cdot)+ \Delta t\big(
      F(v(t,\cdot))+
      \partial_x v(t,\cdot) \mu(t,\cdot)\big)(x)  \Big)_{\R^d}\\
      & \qquad \qquad \forall x\in M, X \in T_xM.
    \end{aligned}
  \end{equation}
  Again, this is a second order elliptic problem of the same type as \eqref{eq:elliptic}, and one can determine consistent initial values for $\mu(0,\cdot)$ from
  \eqref{eq:semidiscrete} without any conflict.
\end{enumerate}
For a comparison, it is instructive to express  the four phase functionals
in terms of the modified vector field $F(v,\mu):=F(v) + \partial_xv \mu$. 
Given $(v,\mu,g)\in \mathcal{U}_d \times \mathcal{V} \times G$ the functionals are defined for
$x\in M,X \in T_xM$ as follows:
\begin{equation*} \label{eq:summarize}
  \begin{aligned}
    \psi_1(v,\mu,g)(x,X) & = (\partial_xv(x)X,F(v,\mu)(x))_{\R^d}, \\
    \psi_2(v,\mu,g)(x,X)& = (\partial_x\hat{v}(x)X,F(v,\mu)(x))_{\R^d},\\
    \psi_3(v,\mu,g)(x,X)& = (\partial_xF(v,\mu)(x)X, \Delta_Mv(x))_{\R^d}
    + (\partial_x v(x) X, \Delta_M F(v,\mu)(x))_{\R^d},\\
    \psi_4(v,\mu,g)(x,X)& = (\partial_x v(x), \Delta_M(v+ \Delta t F(v,\mu))(x))_{\R^d}.
  \end{aligned}
  \end{equation*}




\sect{Applications to flowing closed curves}
\label{sec4}
In this section we consider the simple Example \ref{ex2} and follow
closed curves in a finite dimensional dynamical system \eqref{eq:dynsyst},
i.e.
\begin{equation} \label{edynsyst}
  \dot{u} = f(u), \quad f\in C^{\infty}(\R^d,\R^d).
\end{equation}
Recall the state space  $\mathcal{U}_d=C^{\infty}(S^1, \R^d)$ and the
vector field $F$ from \eqref{e12}.
The Freezing equation \eqref{e37a},  then leads to the following 
hyperbolic PDE system for the functions
$\mu:\mathcal{I} \times S^1 \to \R$ and $v:\mathcal{I} \times S^1 \to \R^d$: 
\begin{align} \label{e44a}
	\partial_t v(t,x) = & f(v(t,x)) + \partial_x v(t,x) \mu(t,x),\quad
	v(0,x)= u_0(x).
\end{align}
The system will be closed by one of the phase conditions
of Section \ref{sec3.2}.
The time dependent diffeomorphism $g:\mathcal{I}\times S^1 \to S^1$
is then determined from the reconstruction equation \eqref{e37b}, i.e. 
\begin{align}
  \label{e44b}
	\partial_t g(t,x) = & \partial_x g(t,x) \mu(t,x), \quad
	g(0,x)= x.
\end{align}
In the following we  describe our numerical setup and then put
the main focus on the choice of  
phase condition, which is crucial for the Freezing method.

\subsection{Relation to the method of characteristics}
Let us note that one may solve the system \eqref{e44a} for a given function
$\mu$ by the method of characteristics (see e.g. \cite[Ch.3.2]{Ev10}). This method
requires to solve for each $x_0 \in M$ the $d+1$-dimensional system
\begin{equation*} \label{characteristic}
  \begin{aligned}
    \xi'(t) = - \mu(t,\xi(t)), & \quad \xi(0)=x_0, \\
    h'(t) = f(h(t)), & \quad h(0)=u_0(x_0).
      \end{aligned}
\end{equation*}
With the solutions denoted by $\xi(t;x_0)$ and $h(t;u_0(x_0))$, one then finds  
the solution $v$ of \eqref{e44a}  from the relation
\begin{equation*} \label{relatechar}
  v(t, \xi(t;x_0))= h(t; u_0(x_0)), \quad t \in \mathcal{I}, x_0 \in M.
\end{equation*}
In our approach, the vector field $\mu$ defining the characteristics
\begin{align*}
  \{(t,\xi(t;x_0))\in \R^{d+1}: t\in \mathcal{I} \}
\end{align*}
is not given {\em a-priori} but determined during
the computation through the phase condition. Thus the orthogonal phase
condition \eqref{eq:orthphase} may be interpreted as selecting characteristics
which force the $v$-profile into an orthogonal direction (see Figure \ref{fig2}). 
Rather than following the solution along a bunch of characteristics we treat
the PDE \eqref{e44a} directly and try to keep the
structure of the underlying manifold $M$ during time evolution.

\subsection{Discretization of the Freezing system}
For our model equation \eqref{e44a} above we write $S^1 = \R / \Z$ and
decompose $S^1= \bigcup_{i=1}^N S_i$ into $N$ segments
$S_i=[x_i, x_{i+1}] \subset [0,1]$ for $i=1,\ldots,N$. 
In the following, indices related to the grid are considered modulo $N$,
which realizes the periodic boundary conditions.
Further, we use equidistant grids:
\begin{equation*}
  \begin{aligned}
    \Delta x & =\frac{1}{N},\quad x_i := (i-1) \Delta x \in [0,1], \quad i = 1,\ldots,N,\\
    t^n & = (n-1) \Delta t, n=1,2,\ldots \quad  \Delta t>0.
  \end{aligned}
\end{equation*}
Our approximate quantities are denoted by
	\begin{equation*}
	v_i^n \approx  v(t^n,x_i),\;  g_i^n \approx g(t^n,x_i), \;  \mu_i^n \approx \mu(t^n,x_i), \; n=1,2,\ldots\; , \; i=1,\ldots,N
	\end{equation*}
with initial data $v_i^0 := u_0(x_i),\; i=1,\ldots,N.$

We solve \eqref{e44a} here by the Lax-Wendroff (LxW) method commonly used for conservation laws with source terms (see \cite{RLV}, \cite{Hudson}). It is convenient to introduce
the first and second order difference quotients
\begin{align*}
  \partial_{\Delta x}v^n_i &=\frac{1}{2\Delta x}(v^n_{i+1}-v^n_{i-1}), \quad
  \partial^2_{\Delta x}v^n_i =\frac{1}{\Delta x^2}(v^n_{i+1}-2 v^n_{i}+v^n_{i-1}),\\
  \partial_{\Delta x}\mu^n_i &=\frac{1}{2\Delta x}(\mu^n_{i+1}-\mu^n_{i-1}), \quad
  \partial_{\Delta t}\mu^n_i =\frac{1}{\Delta t}(\mu^{n+1}_{i}-\mu^n_i)).
\end{align*}
With these settings and the Jacobian $\partial_vf(v)$ of $f$, the fully discrete scheme reads ($i$ runs modulo $N$)
\begin{equation} \label{LxW}
  \begin{aligned}
    v_i^{n+1}& = v_i^n + \Delta t\big[ \partial_{\Delta x}v_i^n
      + f(v_i^n) \big] + \frac{\Delta t^2}{2} \partial_vf(v_i^n)
    \big[f(v_i^n) +2 \partial_{\Delta x}v_i^n \mu_i^n\big] \\
    & + \frac{\Delta t^2}{2} \Big[ \partial_{\Delta x}^2v_i^n (\mu_i^n)^2
     + \partial_{\Delta x}v_i^n \partial_{\Delta x}\mu_i^n \ \mu_i^n +
      \partial_{\Delta x}v_i^n \partial_{\Delta t}\mu_i^n \Big].
  \end{aligned}
\end{equation}
 If $\mu$ is given, the
 LxW scheme is an explicit method which is second
 order consistent in space and time if $f$, $\mu$ and the solution $v$ are sufficiently smooth. Note that \eqref{LxW} collapses to the Taylor method of
 order $2$ for the ODE \eqref{edynsyst} if $\mu \equiv 0$.
 
 This fits well into our setup of $C^\infty$-manifolds and
 $C^{\infty}$-solutions, and a convergence result for the nonlinear LxW scheme
 in case of a constant function $\mu$ is provided in \cite[Ch.4]{Roe19}.
        However, it is also well known that the LxW scheme produces spurious
        oscillations for standard conservation laws when discontinuities, like shock waves, arise; see e.g.\ \cite[Ch 10.3, 17.2]{RLV02}.

        The main departure from the standard case here is that the function $\mu$
        is not given; rather, it is determined by the phase condition.
        Note that the expression $\partial_{\Delta t}\mu_i^n=
        \frac{1}{\Delta t}(\mu_i^{n+1}- \mu_i^n)$ contains the $\mu$-values at
the next time instance at which we want to impose the algebraic constraint.
        
For the coupling we employ an implicit-explicit method (similar to the
IMEX-RK methods in \cite{ARS97}, \cite{Ro19}). First  note that we can
use \eqref{LxW} to express $\partial_{\Delta t}v_i^n$ in terms of
$\partial_{\Delta t}\mu_i^{n}$
as
\begin{equation*} \label{eq:vfrommu}
  \partial_{\Delta t} v_i^n = a_i^n + b_i^n \partial_{\Delta t}\mu_{i}^n,
\end{equation*}
where the terms $a_i^n,b_i^n$ can be read off from \eqref{LxW}.
For the orthogonal phase condition \eqref{eq:orthphase2} we then approximate $\partial_xv(t^n,x_i)$ by the difference quotient
$\partial_{\Delta x}v^n_i$. This leads to the condition
\begin{align*}
  0 = (\partial_{\Delta x}v_i^n, a_i^n + b_i^n \partial_{\Delta t}\mu_{i}^n)_{\R^d},
\end{align*}
from which $\partial_{\Delta t}\mu_i^{n}$ can be determined. With these
values we perform an explicit LxW-step \eqref{LxW}  to obtain $v_i^{n+1}$.

In general, the  phase condition provides the algebraic side constraint for \eqref{LxW} turning it into a
system of $N(d+1)$ equations for the $N(d+1)$ unknowns $v_i^{n+1},\mu_i^{n+1},i=1,\ldots,N$. The IMEX method described above avoids to solve such
 a fully nonlinear system.
 For the energy based phase conditions this implicit-explicit
 coupling is somewhat more involved, see below.
        
        For stability reasons the time-step of the LxW method should be chosen to satisfy
        the standard CFL condition $|\mu| \frac{\Delta t}{\Delta x} < 1$ for fixed $\mu$, see
        \cite[p.110 f.]{RLV}. 
	In our case $\mu$ is the free variable determined by an additional phase condition. 
	Due to this condition it may happen that $\mu$ itself depends on $\Delta t$ causing
	a negative feedback when $\Delta t$ is decreased. 
	Thus in practice we simply choose $\Delta t = \tfrac{1}{2} \Delta x$ and damp $\mu$ uniformly 
	if needed to satisfy the CFL condition.

The reconstruction equation \eqref{e44b} does not involve a source term and is solved in a similar way with less difficulties.
With the half-step values
\begin{equation*} \label{advecRLxWHP}
g^{n+1/2}_{i \pm 1/2} = \frac{g^{n}_{i \pm 1} +
g^{n}_{i}}{2} \pm  \Delta t \frac{g^{n}_{i \pm 1} - g^{n}_i}{2
	\Delta x} \mu^n_i.
\end{equation*}
the LxW scheme for $g$ reads
\begin{equation*} \label{advecRLxW}
g_i^{n+1} = g_i^{n} + \tfrac{\Delta t}{\Delta x} \big(g^{n+1/2}_{i+1/2}-g^{n+1/2}_{i-1/2}\big) \mu^{n+1}_i , \quad n= 1,2,\ldots .
\end{equation*}
Note that for the reconstruction equation we use the boundary
condition
$g(t,1) = 1 + g(t, 0)$ for all $t \in \mathcal{I}$ which lifts the function
$g(t,\cdot)$ from $S^1$ to $\R$.

\section{Numerical Results}\label{sec:results}
\subsection{The normal form system for Hopf bifurcation}
\label{sec4.2}
Our first system is the two-dimensional normal form of the
Andronov-Hopf bifurcation (see e.g.
\cite{Ku04})
\begin{equation} \label{normHopf}
  f(u) = \begin{pmatrix}  - u_2 + u_1(\lambda- (u_1^2+u_2^2)) \\
    \ u_1 + u_2(\lambda-(u_1^2+u_2^2))
    \end{pmatrix}.
\end{equation}


In the following we use $\lambda=1$ so that the unit circle $S^1$  is
an asymptotically stable periodic orbit for \eqref{normHopf}.
As initial data $u_0$ for our flowing curve we choose a circle with center
$(0.6,0)$ and radius $0.6$ (which intersects the system original periodic orbit). 

For comparison, the first simulation is performed using the original system, i.e. we parameterize with time.
Numerically this corresponds to \eqref{LxW} with $\mu \equiv 0$, which
corresponds to the second order Taylor method. Now we compare with the solution of the Freezing system, which uses an arbitrary parameterization by $x \in S^1$, determined during integration.

\begin{figure}
    \centering
    \includegraphics[width=0.9\linewidth]{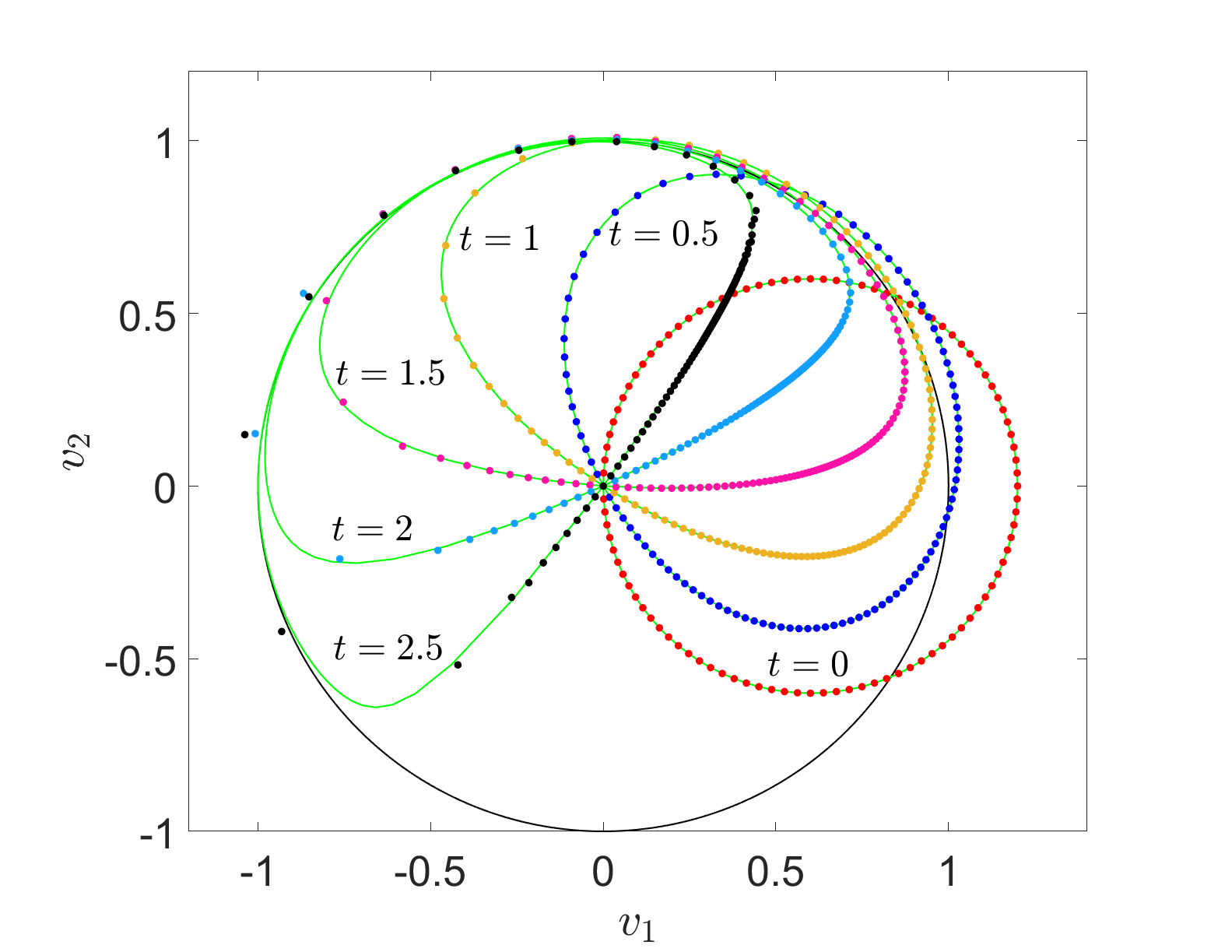}
    \caption{Comparison of flowing curves for the original system
    (continuous lines) at times $t=0,0.5,1,1.5,2$ and $2.5$ with the LxW solution of  
    the Freezing system with orthogonal phase condition (dotted lines showing the distribution of spatial nodes). The initial data $u_0$ for our flowing curve we choose a circle with center $(0.6,0)$ and radius $0.6$.}
    \label{fig5.1}
\end{figure}

In Figure \ref{fig5.1} the continuous curves represent the images of
the initial circle obtained in this way for $t=\frac{n}{2},n=0,\ldots,5$.
We compare them with the images $v^n_i,i=1,\ldots,100$ obtained by solving
the Freezing system with time-step $\Delta t=10^{-3}$ invoking the
orthogonal phase condition. Note that this time is different from the time used in the original system, $\dot{u}=f(u)$, with $f(u)$ given by \eqref{normHopf}. 

While the overlap
is almost perfect for the first time instances, there are slight deviations
at later times $t=2,2.5$, in particular at the extremal points to
the left. 
In this region, neighboring image points are moved apart as
shown by the distribution of nodes provided by the Freezing approach,
so that linear interpolation of neighboring nodes becomes problematic for
graphical representation.



\begin{figure}    
    \begin{tabular}{cc}
     \includegraphics[width=0.5\linewidth]{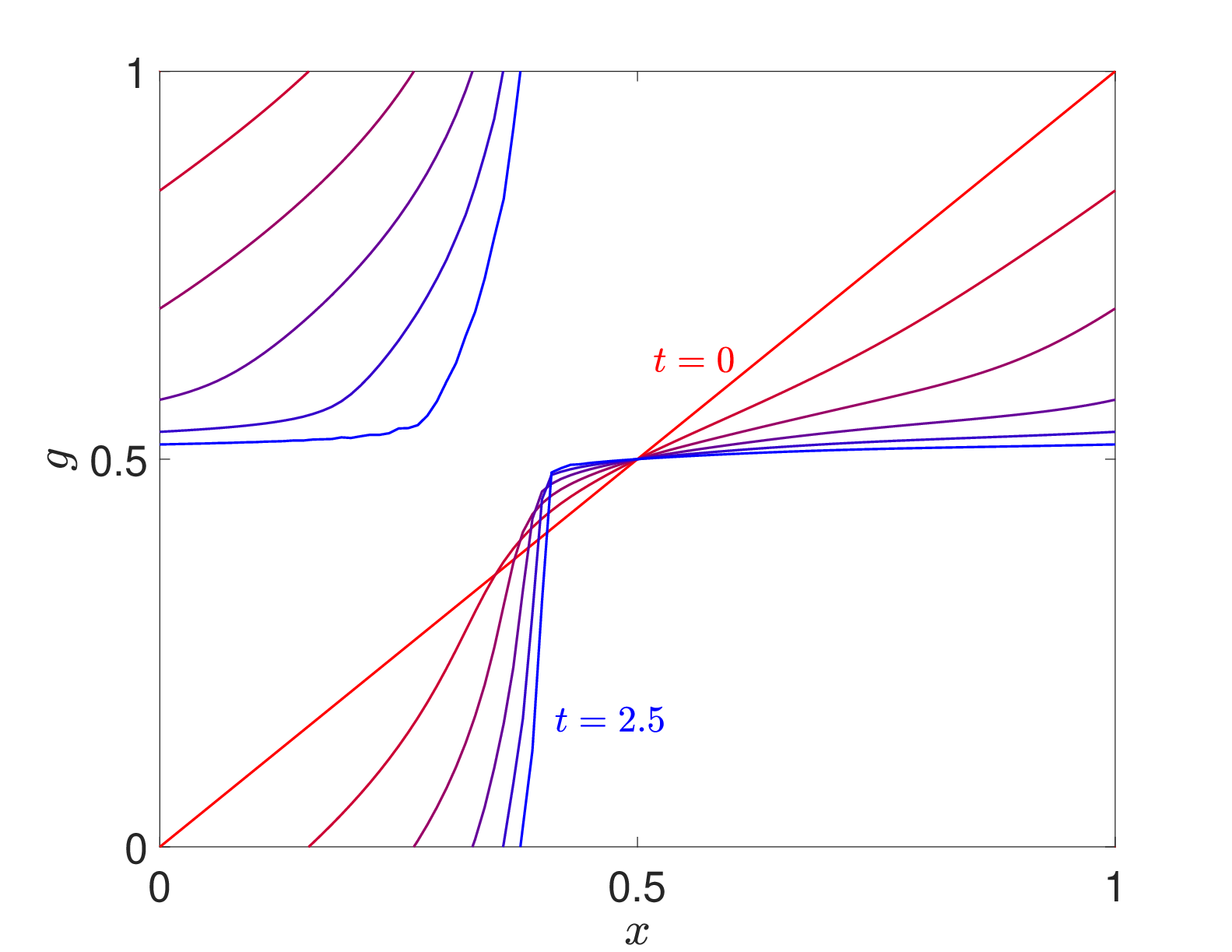}    &  \includegraphics[width=0.5\linewidth]{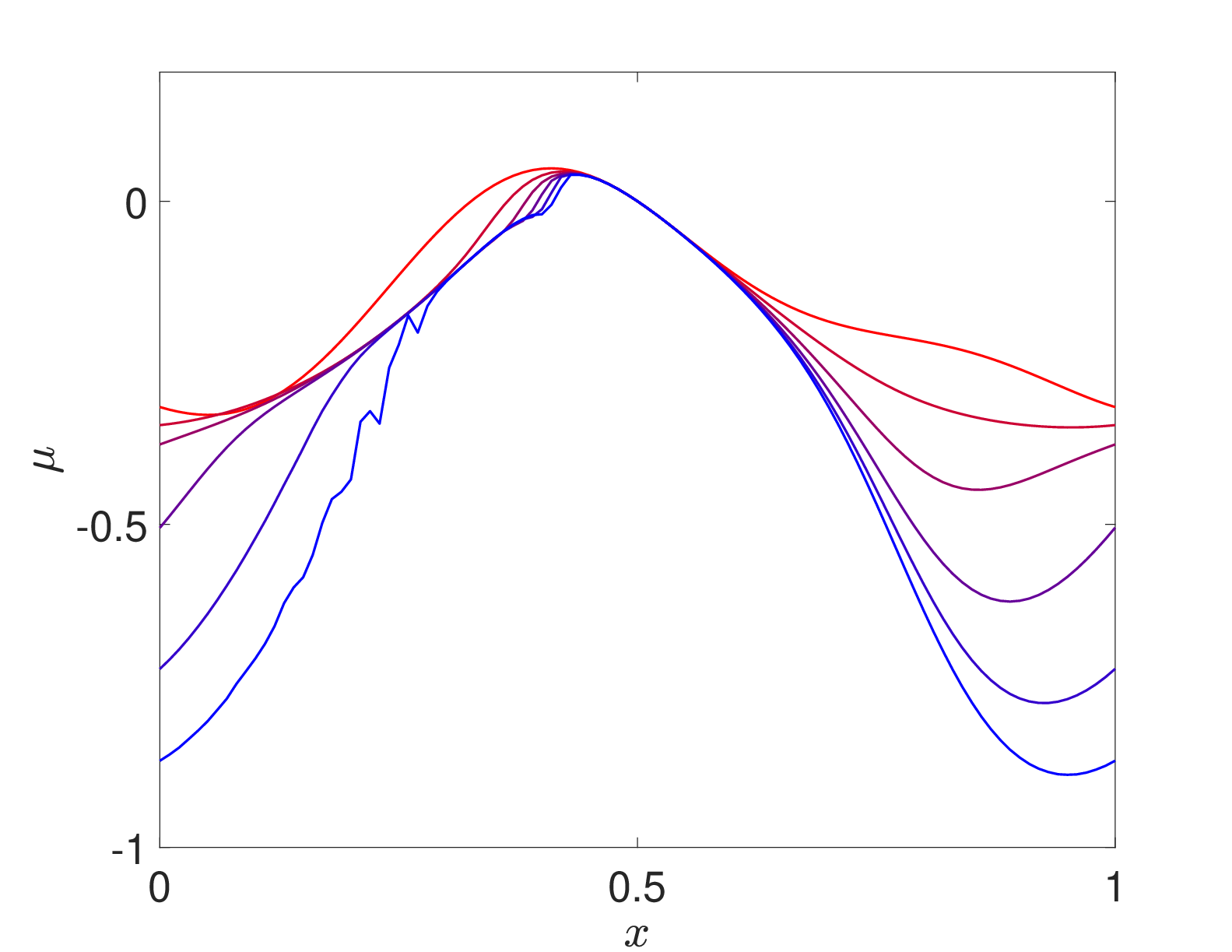}\\
     (a)    & (b) 
    \end{tabular}
    \caption{The left graph shows the solutions $g$ of the reconstruction equation for the orthogonal phase condition for the times $t=0$ (red curve)
    and successively for $t=0.5,1,1.5,2$ and $2.5$. The corresponding curves for the
    control variables $\mu$ are shown in the right graph.
    Note that the graph of the diffeomorphism becomes critical in very flat
    and very steep regions, which occur in the mid-interval as
    time progresses.}  
    \label{fig5.2}
\end{figure}

For the same experiment we show in Figure \ref{fig5.2}  the functions $\mu(t,\cdot)$ and the reconstructed diffeomorphism $g(t,\cdot)$ as obtained from
the reconstruction equation.

Although $\mu$ appears smooth, the reconstructed diffeomorphism $g$ for this simulation exhibits intervals 
with critical steep regions ($\partial_x g \rightarrow \infty$) as well as flat regions ($\partial_x g \rightarrow 0$)
which tend to violate the invertibility of $g$. 
This points towards the use of phase conditions which also involve gradients of the solution;
see Section \ref{sec4.4}.

In a second experiment we follow an initial circle which intersects the
periodic orbit and which is deformed over time into a  stretched banana-like
curve, shown as a $3D$ graphic in Figure \ref{fig1}.
The first graphs show that the deviation of the exact continuous curves
from the numerical solution of the Freezing system with orthogonal phase conditions
starts being visible at even earlier times than in Figure \ref{fig5.1}.

\begin{figure}[!ht]\label{fig5.3}
    \begin{center}
      \includegraphics[width=0.9\textwidth]{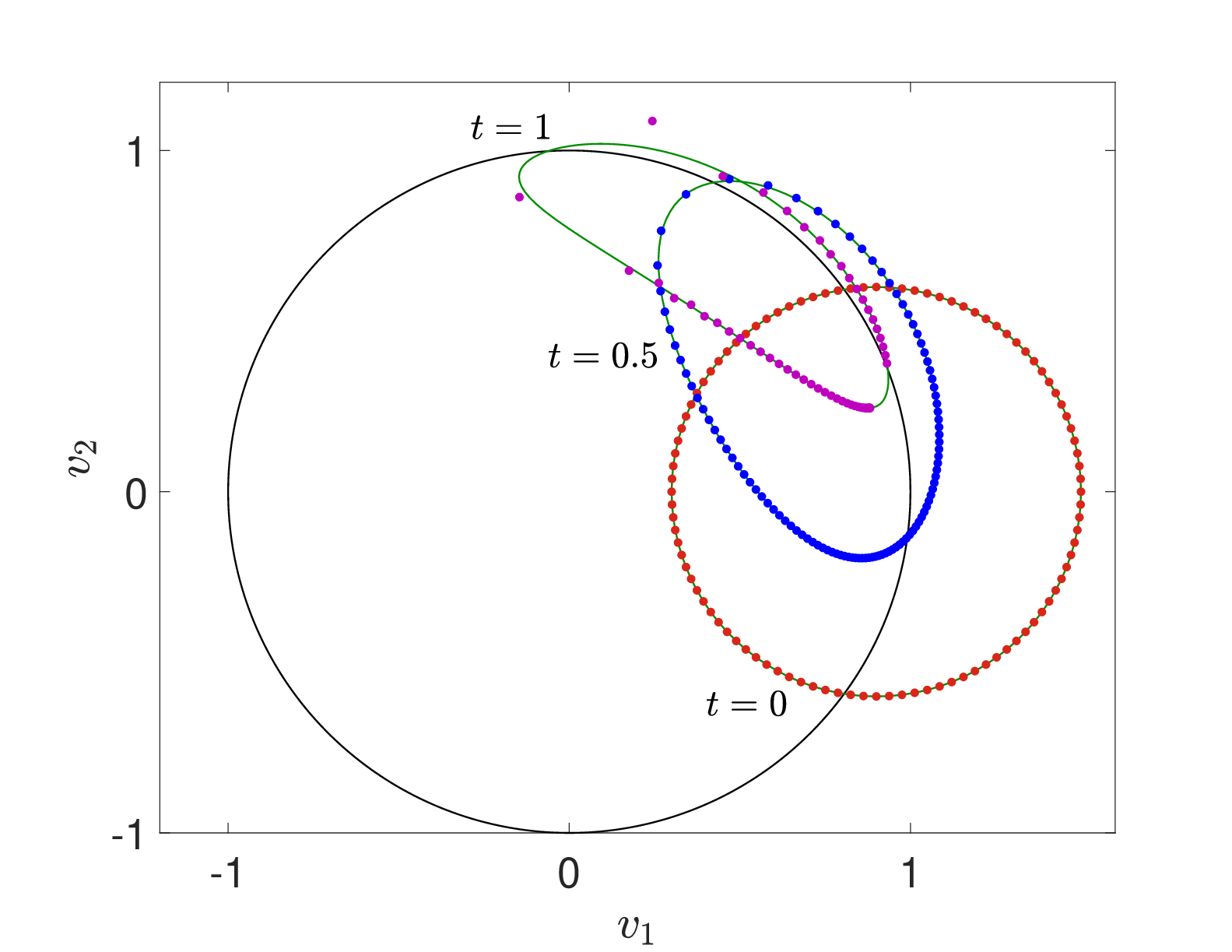}
    \end{center}
    \caption{Comparison of flowing curves for the original system
    (continuous lines) at times $t=0,0.5$ and $1$ with the LxW solution of  
    the Freezing system with orthogonal phase condition (dotted lines showing the distribution of spatial nodes). The initial data $u_0$ for this flowing curve is a circle with center $(0.9,0)$ and radius $0.6$.}
\end{figure}

\subsection{Energy based phase conditions}
\label{sec4.4}
\begin{figure}
    \centering
    \includegraphics[width=0.9\linewidth]{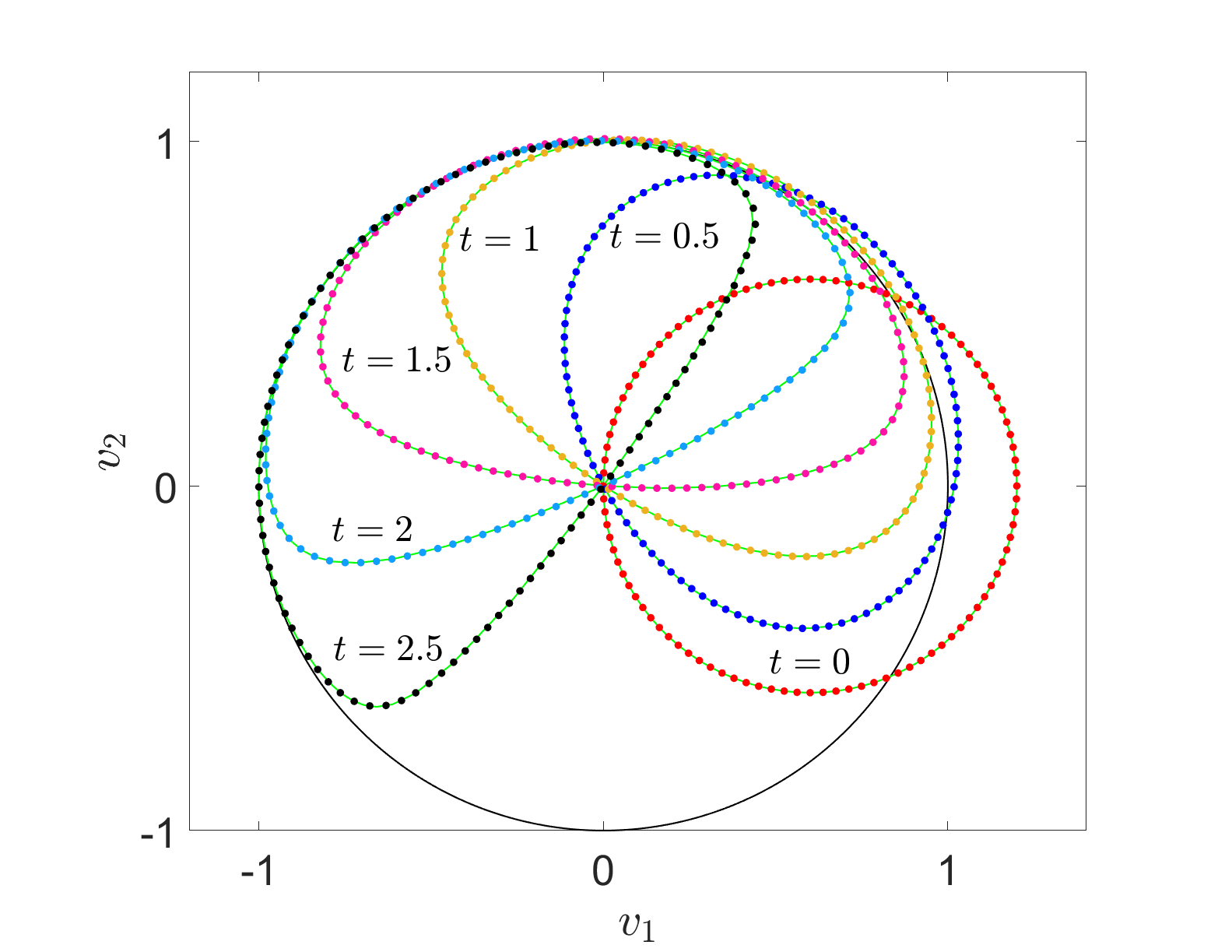}
    \caption{Comparison of flowing curves for the original system
    (continuous lines) at times $t=0,0.5,1, 1.5, 2$ and $2.5$ with the LxW solution of  
    the Freezing system with semi-discrete energy based phase condition (dotted lines showing the distribution of spatial nodes).}
    \label{fig5.4}
\end{figure}

While the LxW scheme \eqref{LxW} is retained, we determine the $\mu$-values
at time $t^{n+1}$ from the semi-discrete condition \eqref{eq:semidiscrete}.
In fact, equation \eqref{eq:semidiscrete} is written down for all $x_i,
i=1,\ldots,N$ at $t^{n+1}$ with all difference quotients of $v$ evaluated
at time $t^n$.
Then one obtains a linear system for the values $\mu_i^{n+1},i=1,\ldots,N$ which
is a discretization of the two-point boundary value problem \eqref{eq:semidiscrete} with periodic boundary conditions. The associated matrix is almost
tridiagonal with some extra entries in the upper right and lower left
corners. Pursuing this approach for the first example in Section \ref{sec4.2} (see Fig. \ref{fig5.1})
yields visually perfect coincidence with the exact continuous curves; see Figure
\ref{fig5.4}.

The following Figure \ref{fig5.5} shows the development of the diffeomorphism $g(t,\cdot)$
and the control variable $\mu(t,\cdot)$:

\begin{figure}
    \begin{tabular}{cc}
      \includegraphics[width=0.5\linewidth]{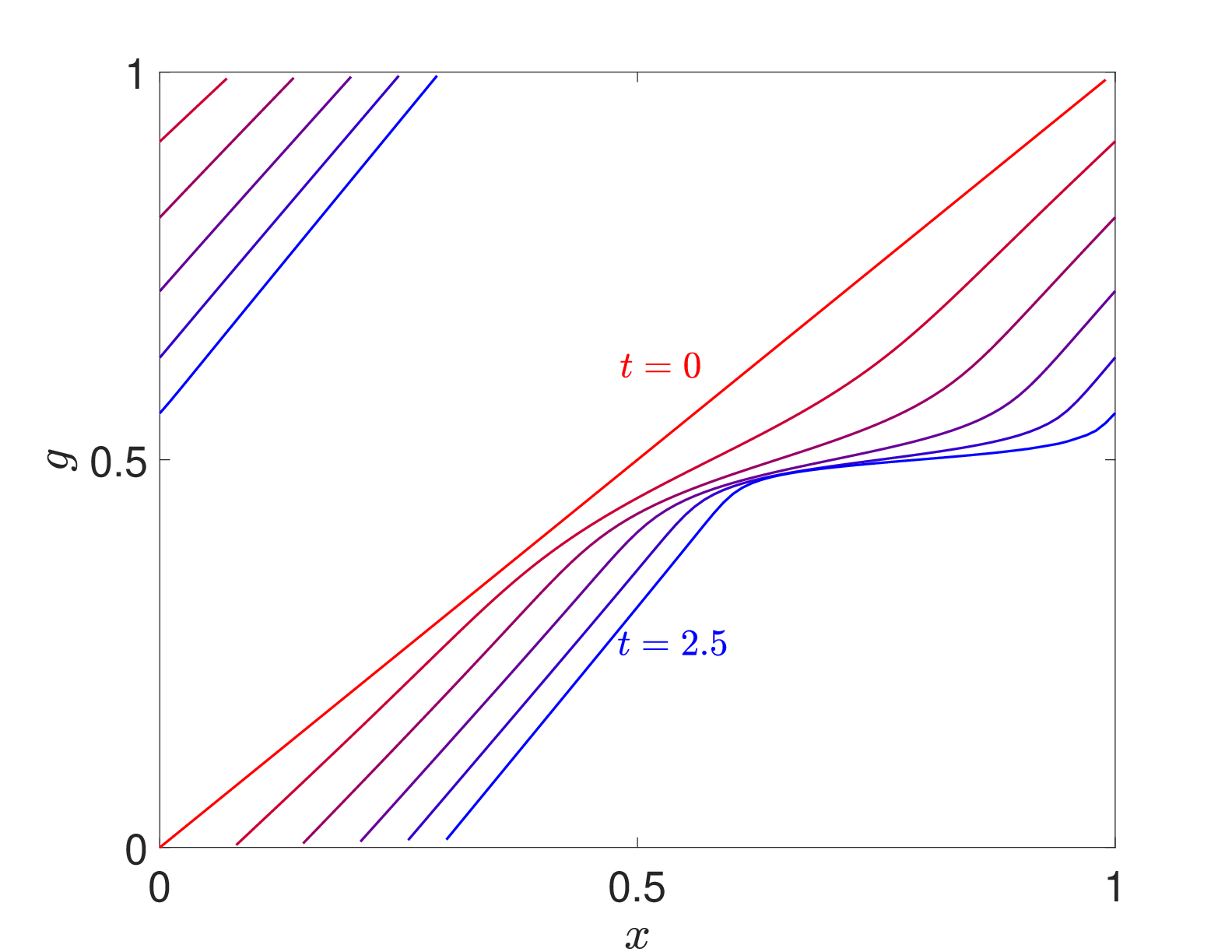}   &  \includegraphics[width=0.5\linewidth]{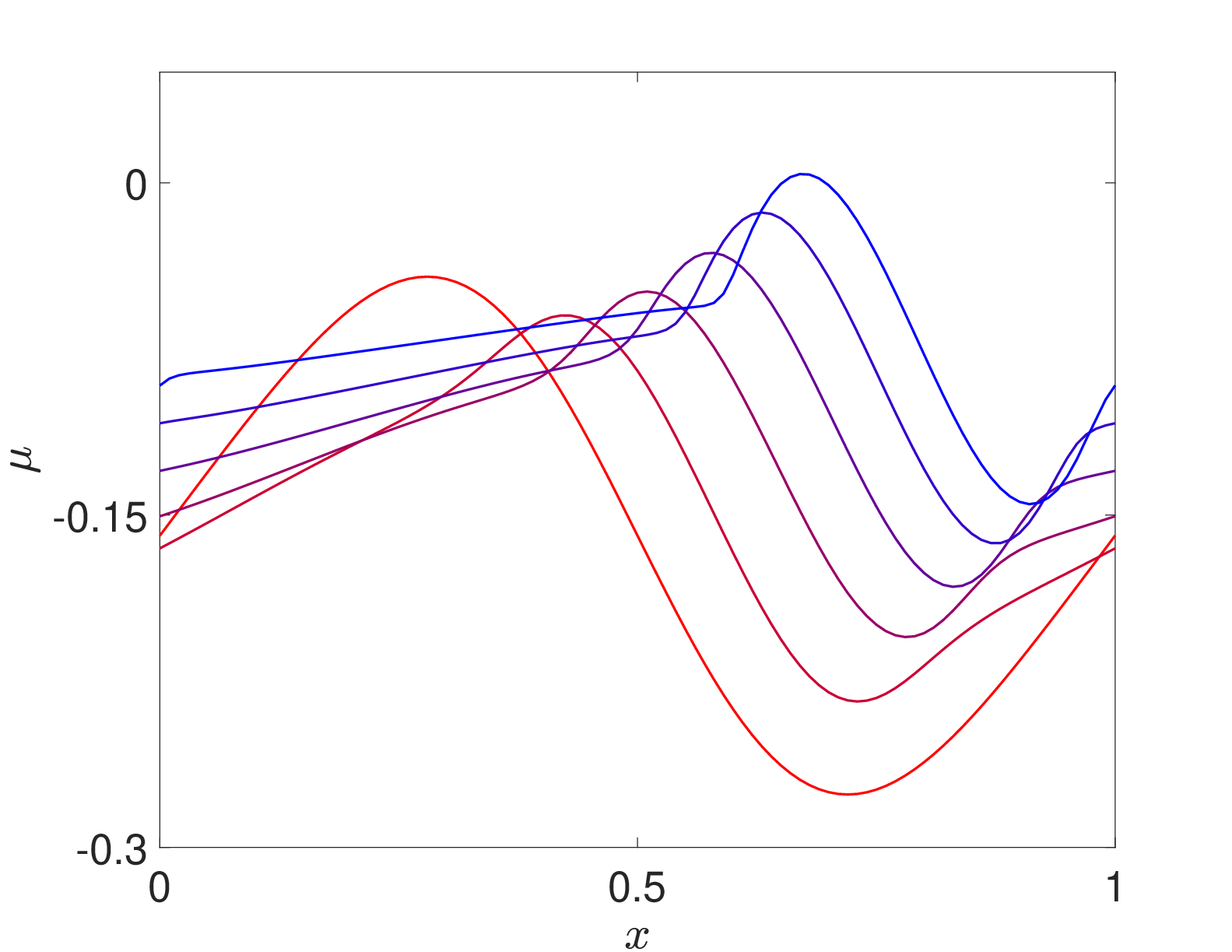} \\
      (a)   &  (b)
    \end{tabular}
    
    \caption{Evolution of the diffeomorphism $g(t,\cdot)$ for the LxW scheme
      with semi-discrete energy based phase condition at times
      $t=0,0.5,1,1.5,2$ and $2.5$ for the solutions shown in Figure \ref{fig5.4} (left)
      and evolution of the control variable $\mu$ for the same data (right).}
    \label{fig5.5}
\end{figure}


\begin{figure}[!ht]
    \begin{center}
      \includegraphics[width=0.9\textwidth]{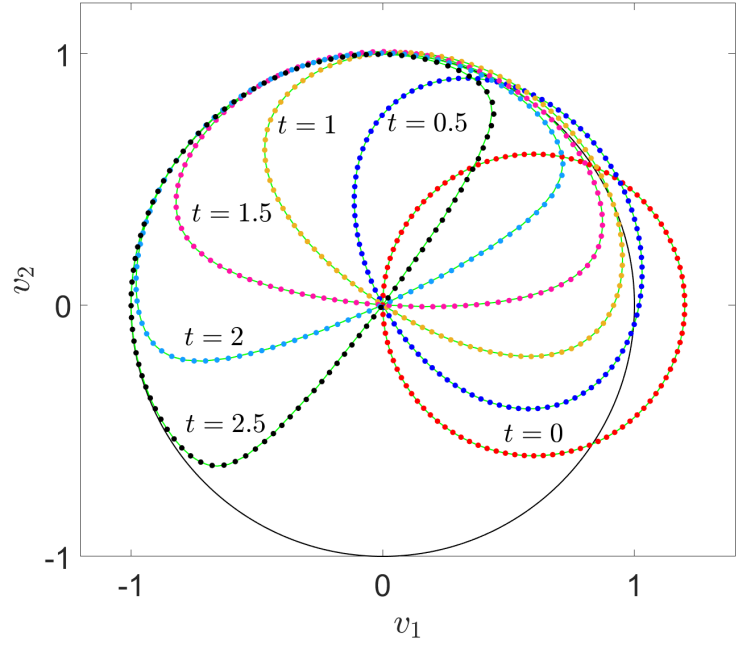}
    \end{center}
    \caption{Comparison of flowing curves for the original system
    (continuous lines) at times $t=0,0.5,1, 1.5, 2$ and $2.5$ with the LxW solution of  
    the Freezing system with semi-discrete energy based phase
    condition (dotted lines showing the distribution of spatial nodes). }
    \label{fig5.6}
\end{figure}

The superiority of the energy based phase condition over the
orthogonal phase condition is also clearly
seen for the second example from Figure \ref{fig5.3} where the initial
curve is strongly distorted; see Figure \ref{fig5.6}. We also note
the almost uniform distribution of nodes on the curve. 

\begin{figure}
    \begin{tabular}{cc}
     \includegraphics[width=0.5\linewidth]{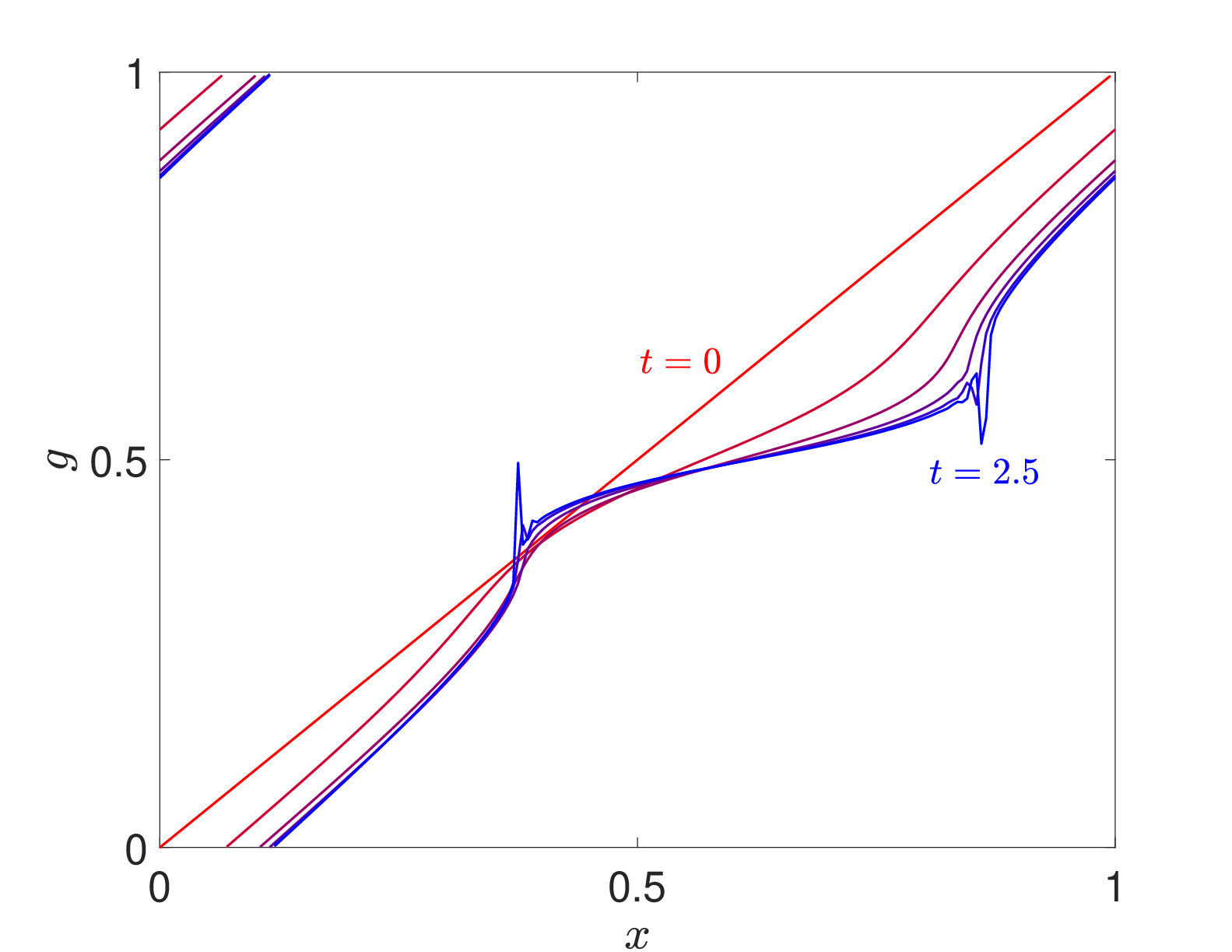}    &  \includegraphics[width=0.5\linewidth]{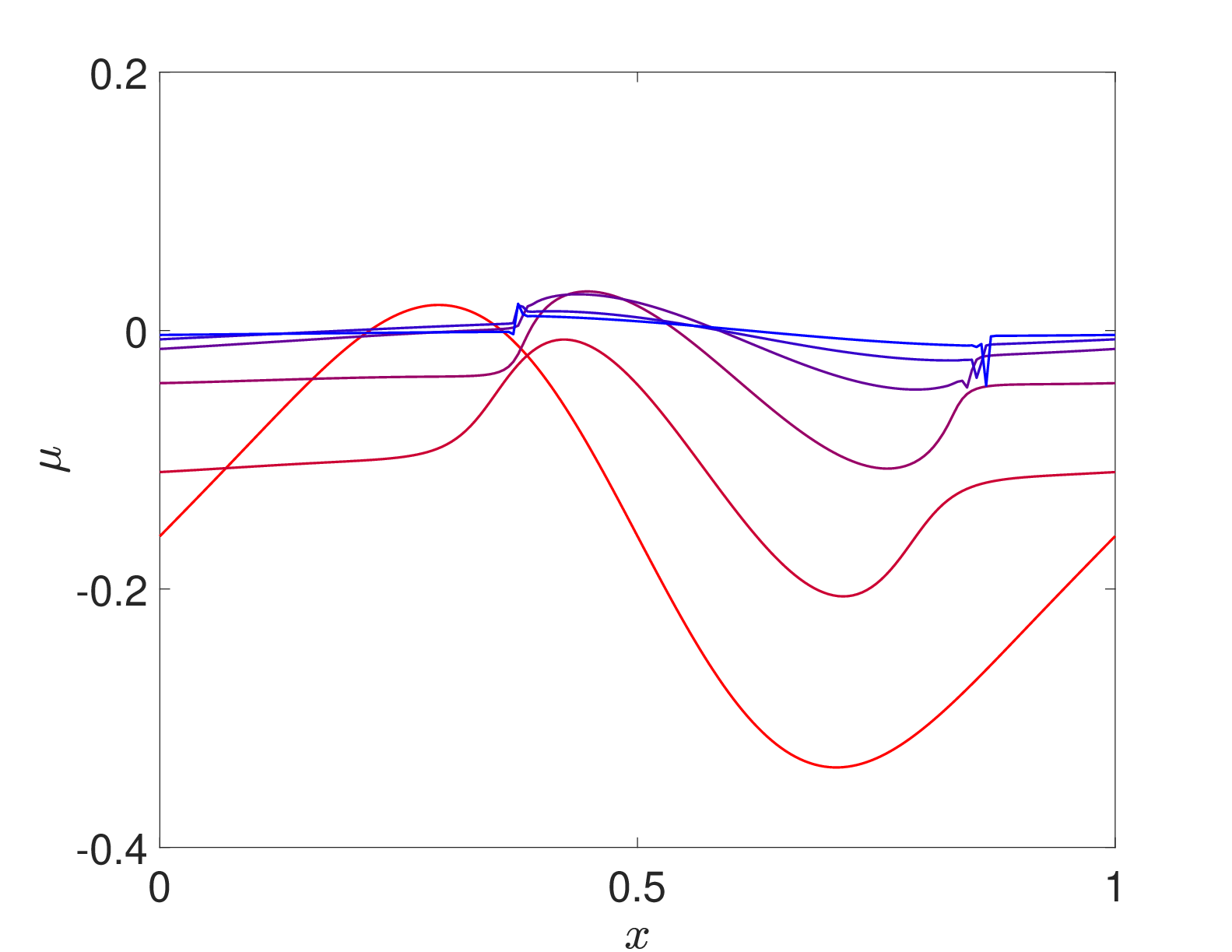} \\
      (a)   & (b)
    \end{tabular}
    
    \caption{Evolution of the map $g(t,\cdot)$ for the LxW scheme
      with semi-discrete energy based phase condition at times
      $t=0,0.5,1,1.5,2$ and $2.5$ for the solutions shown in Figure \ref{fig5.6} (left)
      and evolution of the control variable $\mu$ for the same data (right). }
    \label{fig5.7}
\end{figure}

The corresponding evolution of the map $g(t,\cdot)$ and the control variable
$\mu(t,\cdot)$ is shown in Figure \ref{fig5.7}.
Note that the diffeomorphic property is finally lost at times when the
curve develops very sharp bends (high curvature).
This problem can be alleviated to a certain extent by passing from the
energy integral to a curvature integral. We do not
discuss the detailed results here but refer to \cite[Ch.5.2.4]{Roe19}.

\subsection{Grid evolution on energy levels}
\label{sec4.3}
Our next example is a Hamiltonian system \eqref{edynsyst}  with the vector
field of  the cubic oscillator 
	\begin{equation} \label{cubic} 
	f_{\mathrm{cubic}} (u) = \begin{pmatrix} u_2 \\ -4 u_1^3 + 2 u_1
    \end{pmatrix} .
	\end{equation}
Level curves of the energy 
\begin{equation*} \label{CubicEnergy}
E(u) = \frac{1}{2} u_1^2 +V(u_2), \quad V(u_2) = u_2^2(u_2^2-1).
\end{equation*}
are invariant under the flow and hence fixed points of the evolution operator \eqref{e12}. 
This allows us to compare the distribution of nodes for the various phase conditions.
First  note that the level curves are well approximated by our LxW method \eqref{LxW}
due to its second order accuracy.

The initial data $u_0$ are given  by $u_0 (x_i) := v_{\mathrm{sol}} (i \Delta \tau)$ where
$v_{\mathrm{sol}}$ is obtained from a standard ODE integrator 
with step-size $\Delta \tau$ applied to \eqref{edynsyst}, \eqref{cubic}.
Since this deviates from arclength parameterization, the  corresponding nodes are
non-uniformly distributed on the image curve.

\begin{figure}
     \begin{tabular}{cc}
     \includegraphics[width=0.5\linewidth]{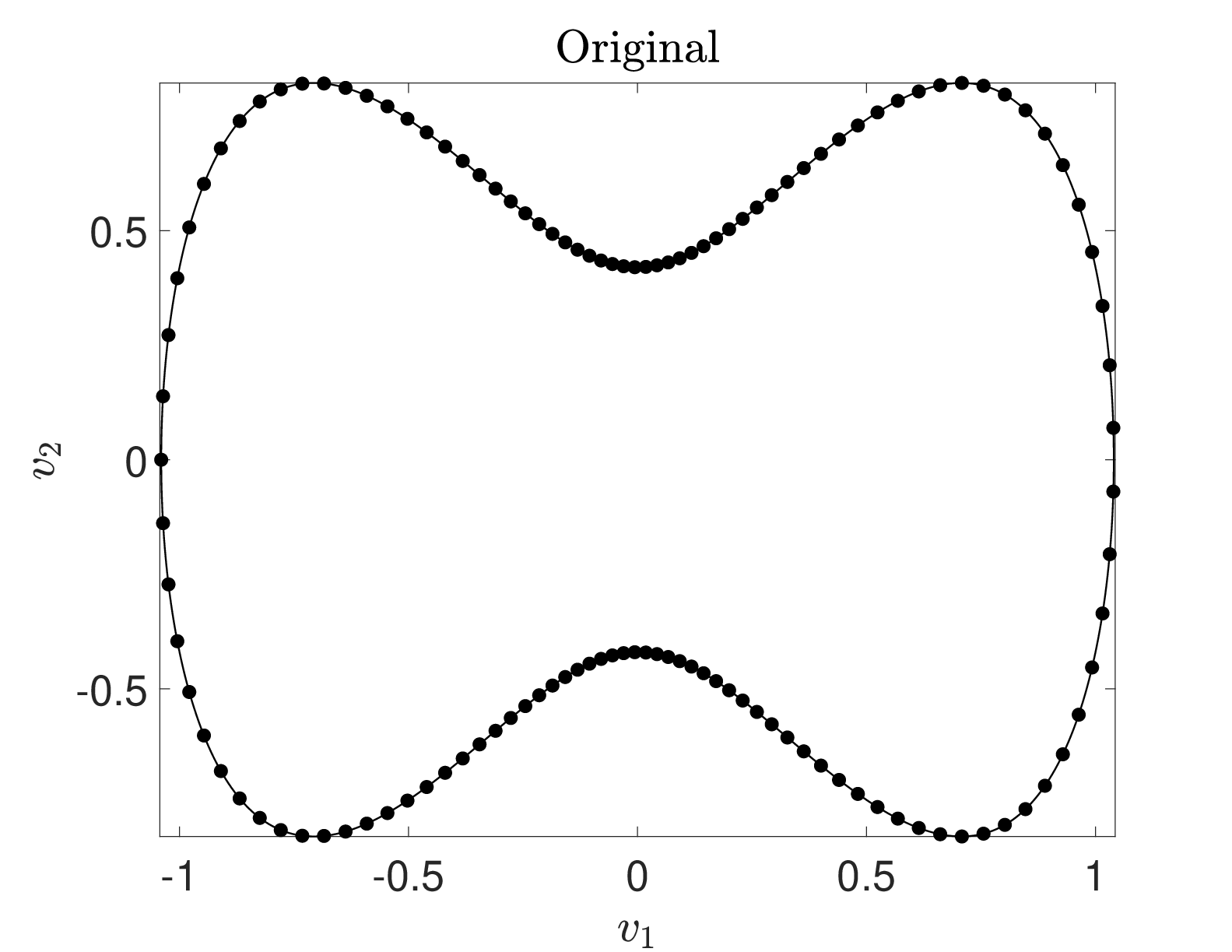}    &  \includegraphics[width=0.5\linewidth]{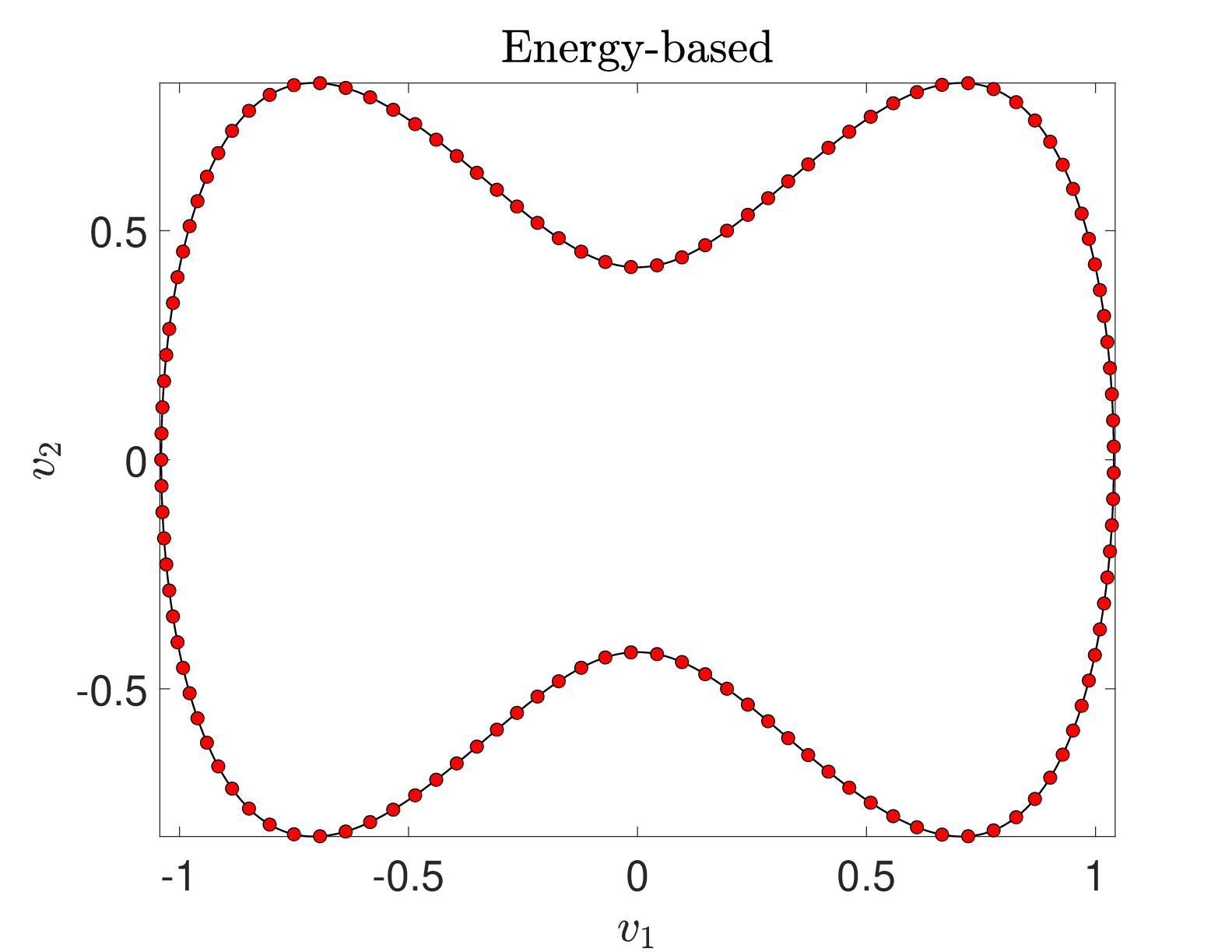} \\
      (a)   & (b)
    \end{tabular}
    \caption{Two nodal distributions on an energy level of the cubic oscillator.
    The left panel emerges from a fixed time step integration of
    original system ($\mu \equiv 0$) while the right panel distribution emerges from
    solving  the Freezing system with semi-discrete energy based phase condition.}
    \label{fig5.8}
\end{figure}

\begin{figure}[!ht]
    \begin{tabular}{cc}
    \includegraphics[width=0.5\linewidth]{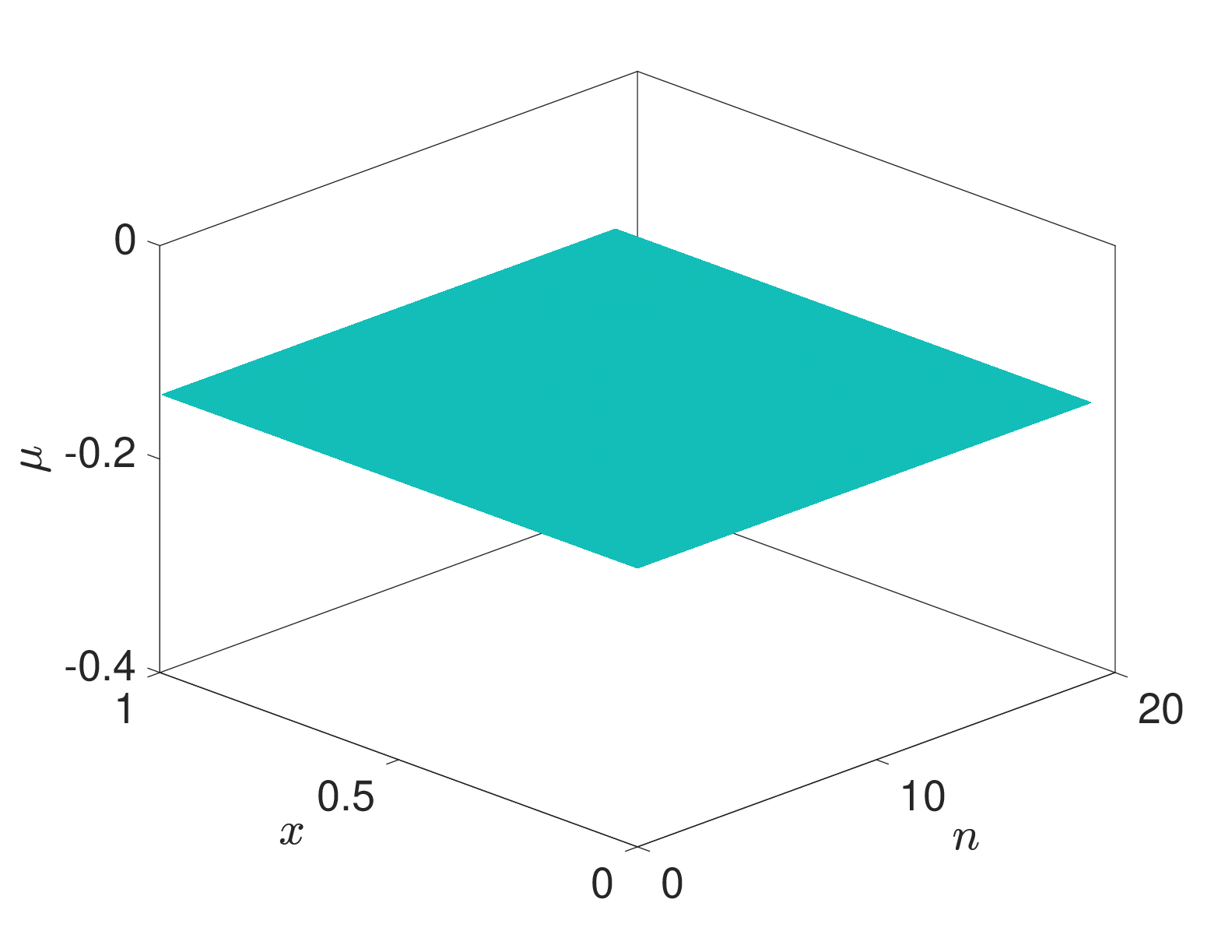}     &  \includegraphics[width=0.5\linewidth]{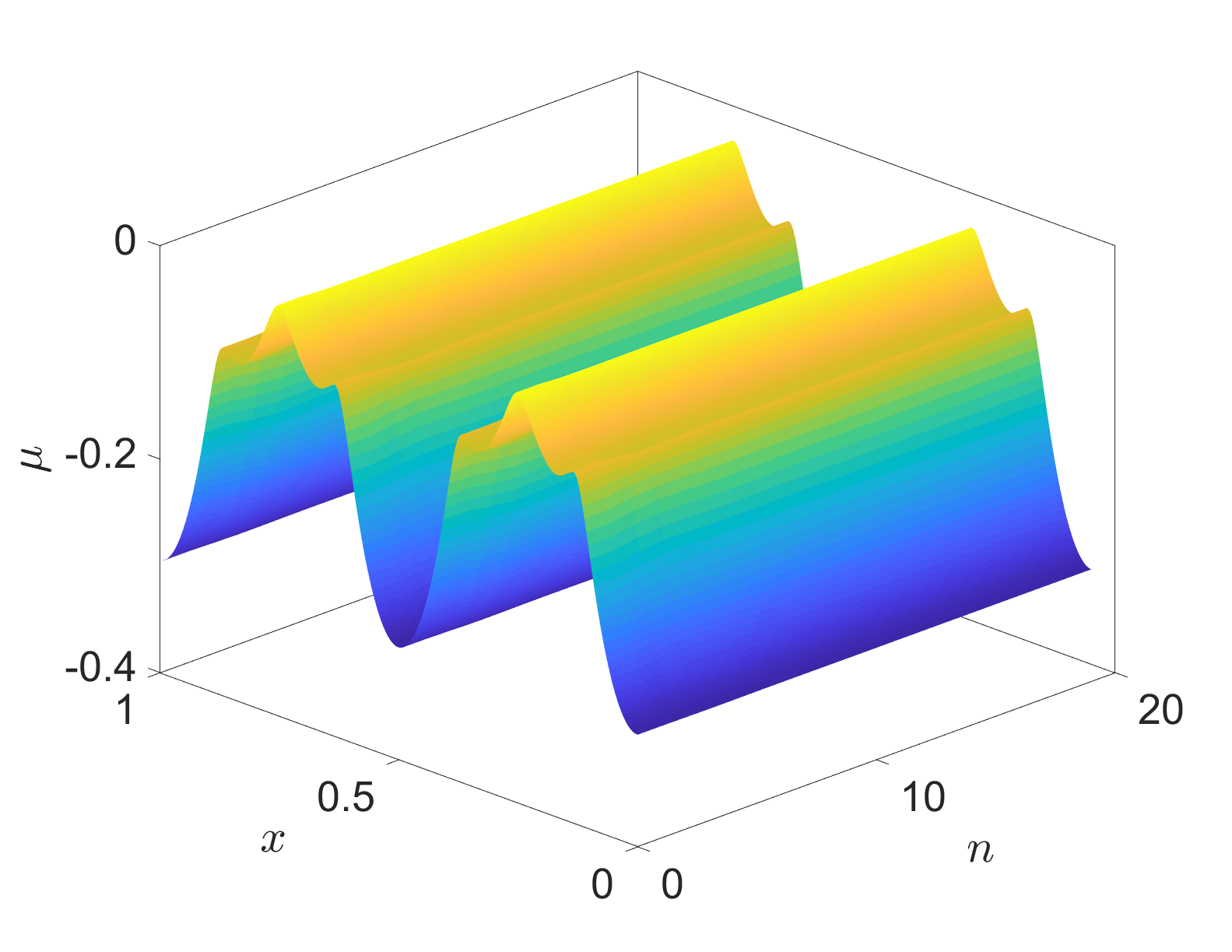} \\
    (a)     &  (b)
    \end{tabular}    
    \caption{Cubic oscillator: Evolution of the $\mu$-values for the original system by equidistant integration in time (left panel) and by the Freezing methods with semi-discrete energy based phase condition (right panel).}
    \label{fig5.9}
\end{figure}

Due to the construction of the initial data this parameterisation is retained naturally as we evolve the original system (this holds exactly  if the step-size in time is
an integer fraction of the period).
On the other hand, the semi-discrete energy based phase condition quickly
adapts the parameterisation
and achieves a uniform distribution of nodes in the image space.
This is confirmed by the corresponding time dependence of the nodes and of the control variable $\mu$ shown in Figure \ref{fig5.9}.
Finally, we note that the orthogonal phase condition leaves all nodes fixed
in this particular example since  
$f$ and $v_x$ are linearly dependent. Thus there is no orthogonal part (see Fig. \ref{fig2}) and 
$\mu$ can be chosen here to nullify the flow field completely (to actually ``freeze" it).

The next Figure \ref{fig5.10} shows the evolution of the difference
quotient $\| \partial_{\Delta x}v^n_i  \|$ for the semi-discrete
energy based phase condition. It shows the quick process which evens out
the spatial distances.
Finally, we show in Figure \ref{fig5.11} the evolution of the diffeomorphism $g$ and the control
variable $\mu$ when solving the Freezing system with the energy based
phase condition. One observes that the $\mu$-variable moves like a traveling
wave while the dynamics is rotating around the level curve.
This time $g(t,\cdot)$ stays smooth and invertible throughout the simulation.
\begin{figure}[!ht]
    \centering
    \includegraphics[width=0.6\linewidth]{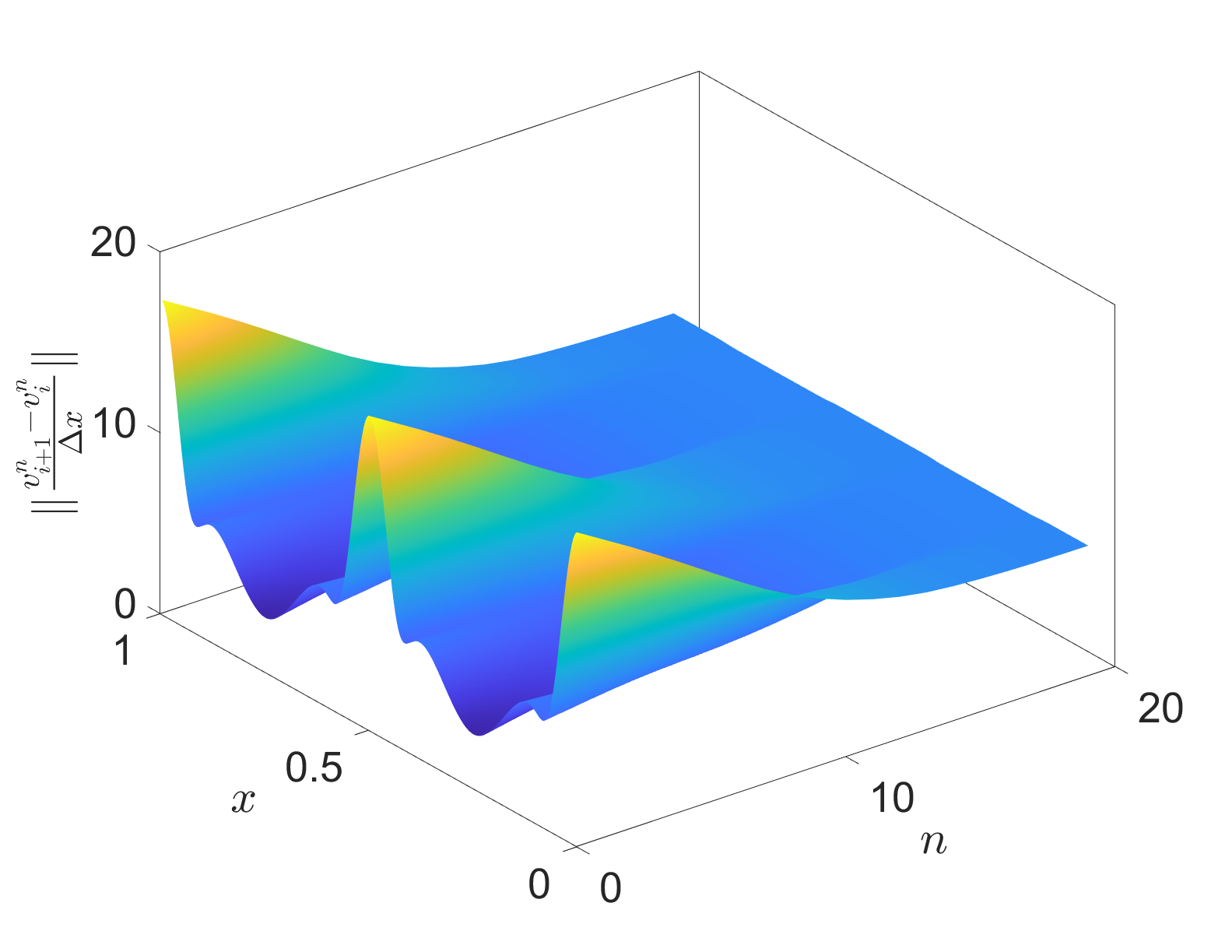}
    \caption{Evolution of difference quotients $||\partial _{\Delta x} v_i^n ||$ during the first 20 steps for the semi-discrete energy based phase condition.}
    \label{fig5.10}
\end{figure}

\begin{figure}[!ht]
\begin{tabular}{cc}
\includegraphics[width=0.5\linewidth]{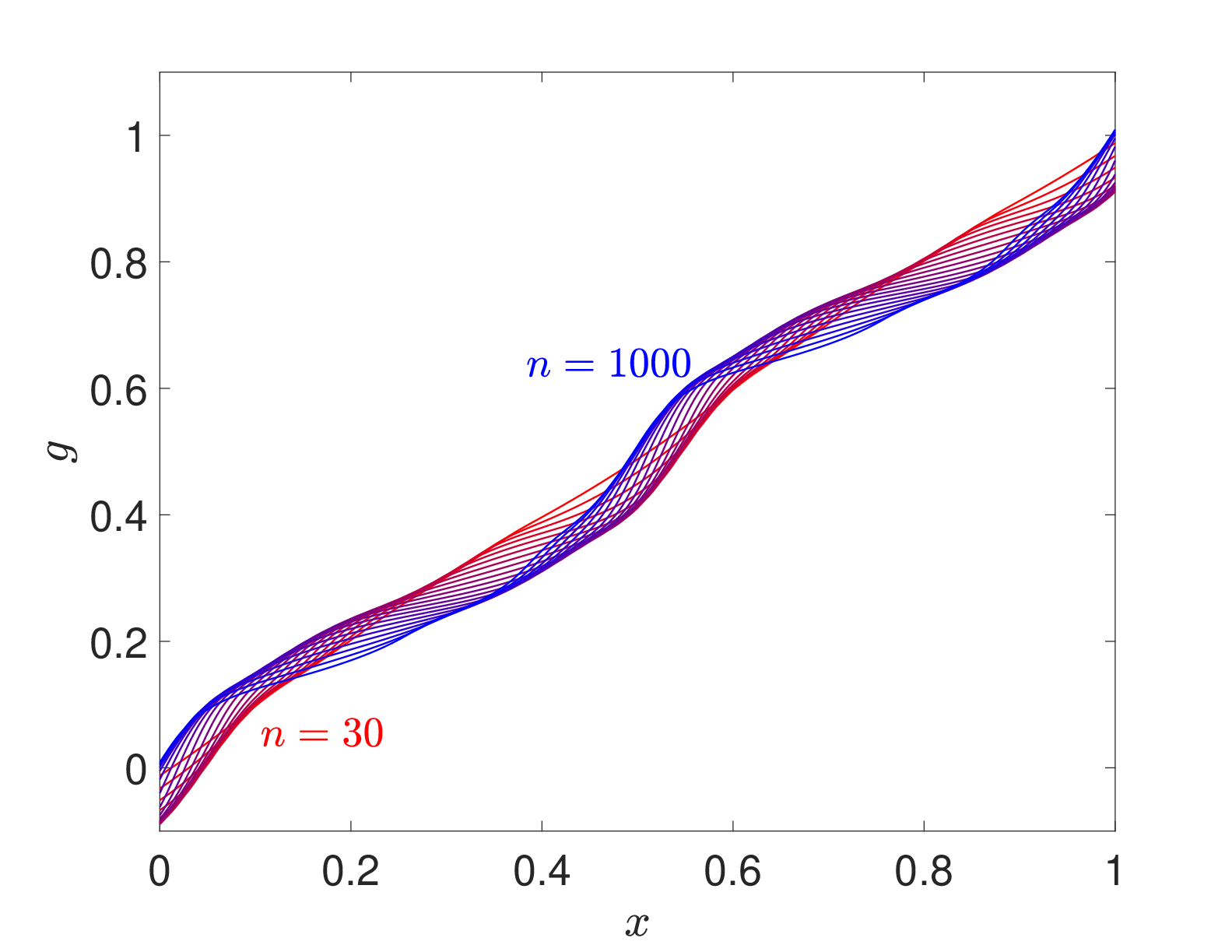}     &  \includegraphics[width=0.5\linewidth]{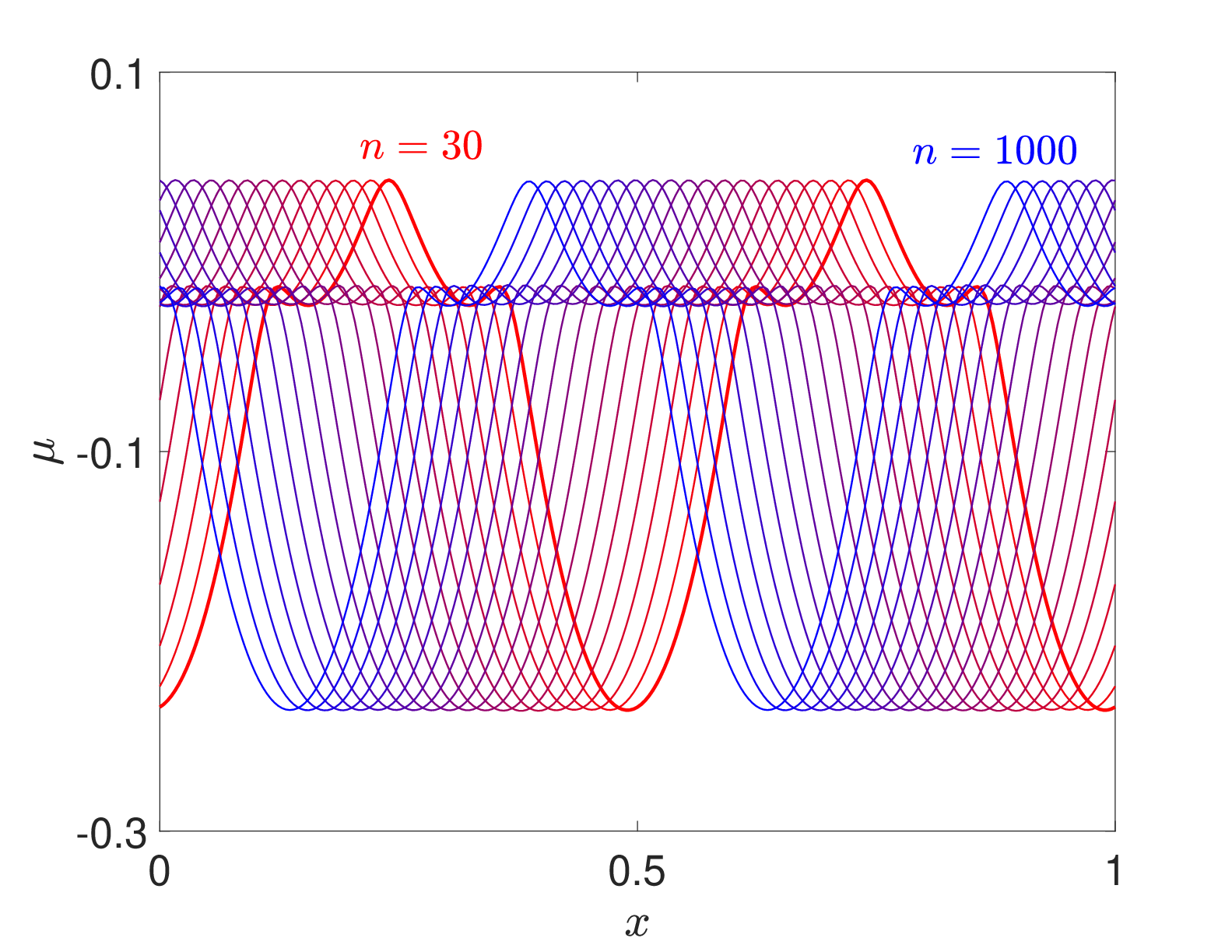}\\
(a)     & (b)
\end{tabular}
    \caption{(a) Evolution of the spatial profile of the diffeomorphism $g$ for the Freezing system with
  energy based phase condition. (b) Evolution of the spatial profile of $\mu$ for the Freezing system with
  energy based phase condition. }
    \label{fig5.11}
\end{figure}

%


\sect{Conclusion and further applications}
\label{sec6}
This paper deals with first order evolution equations
for vector-valued functions defined on a manifold. The vector field
of the evolution equation is assumed to be equivariant under spatial diffeomorphisms of the manifold, i.e. under  arbitrary smooth coordinate transformations.
Our numerical approach exploits this structure by
 introducing new time-dependent variables:  a transformed solution
and a time-dependent diffeomorphism that starts from the identity and is 
driven by an initially unknown time-dependent vector field.  The resulting
extended evolution system \eqref{e1b}, \eqref{e1c} is closed by
imposing side constraints (phase conditions) that effectively determine the additional vector field. In Section \ref{sec2}
we derive several  phase conditions from minimization
principles. Their purpose is to minimize the time derivative of the
solution, making it as stationary as possible, or to minimize
energy-like integrals of the solution, thereby reducing its spatial variations.
In summary, the resulting system to be solved is a PDAE, that is, a partial differential equation with an infinite number (a discretized field) of algebraic side constraints.

We formulate our approach within the framework of evolution 
equations that are equivariant under diffeomorphisms, using the tools
of global analysis for  infinite dimensional Lie groups \cite{KM97}.
From this perspective, our approach can be considered
as an infinite-dimensional generalization of methods derived for
finite dimensional symmetries. We remind the reader that, in the literature, these are known 
as reduction and reconstruction \cite{RKML03, kavousanakis2025going},
Freezing \cite{BT04}, or method of slices \cite{FC12}.

In this paper we restrict our attention  to the simplest setting,  namely the evolution of a parameterized
 manifold $u_0:M\to \R^d$ in the phase
space of a finite dimensional dynamical system
\begin{equation*} \label{eq:dynfinal}
\frac{\partial u}{\partial t}(t,x)= f(u(t,x)), \quad u(0,x)=u_0(x), x\in M, \quad f \in C^{\infty}(\R^d,\R^d).
\end{equation*}
In this case the PDAE becomes a symmetric hyperbolic  system
\eqref{e44a}, \eqref{e44b} with side constraints. We solve it
by  the explicit Lax-Wendroff method supplemented
by  an implicit step that enforces  the algebraic constraints.
We apply the method to closed curves in the phase space
of two types of oscillators (limit cycles).  The results show that energy-based
phase conditions offer a clear advantage over time-based conditions,
when the curves in phase space do not converge but instead become 
strongly distorted.

We view these simulations as test cases for 
more general systems in which the vector field in  function space depends on spatial derivatives, in particular on second-order
derivatives. In the following, we mention several examples which fit
into the abstract framework and to which our method applies.

The first example is the curve shortening equation
\begin{equation} \label{cs}
u_t = \frac{1}{|u_x|^2}P(u_x)u_{xx} =:F(u),
\end{equation}
where $u:M\to \R^d$ with $M=S^1$ or $M=\R$ and $P(v)=I_d -
\frac{1}{|v|^2}v v^{\top}$ is the orthogonal projector onto
$v^{\perp}$, $v \in \R^d$; see \cite{AAAW13}.
Note that \eqref{cs} is a coordinate free version of the
implicit form given in \cite{AAAW13} 
\begin{equation} \label{csal}
P(c_x)c_t= c_{xx}, 
\end{equation}
which parameterizes
the curve $c:M \to \R^d$ by arclength.
One readily verifies that  $F$ in \eqref{cs} is equivariant
under the transformation from $u(t,x)$ to $u(t,g(t,x))$ for $g(t,\cdot)
\in \mathrm{Diff}(M)$, and that the arclength version \eqref{csal} is obtained by choosing $g(t,x)=\int_0^x |u_x(t,\xi)| \mathrm{d \xi}$.
We emphasize that the Freezing system \eqref{e37a}, \eqref{e37b}
for the curve shortening equation remains parabolic,
 since its principal part is elliptic and the additional
terms are only first order in space.

Similar assertions apply to the second order equation which
describes the motion of a vortex filament $u(t,\cdot):M \to \R^3$ ; see \cite{Kh13}
\begin{equation*} \label{eq:vf}
u_t= \frac{1}{|u_x|^3} u_x \times u_{xx}=: F(u).
\end{equation*}
Again, this is the coordinate free version of the equation
$P(c_x)c_t = c_x \times c_{xx}$ used in \cite{Kh13} for curves parameterized
by arclength.

 As a final application we mention mean curvature flow. To be  more specific,
 consider the motion of a   $p$-dimensional submanifold of $\R^d$
 parameterized by some  $C^{\infty}$-map $u(t,\cdot): M \to \R^d$ on an open
 subset $M \subseteq \R^p$. The equation is
 \begin{equation} \label{eq:mcf}
 u_t= \Delta_M u=:F(u) ,\quad u(0,x)=u_0(x), \quad x\in M,
 \end{equation}
where $\Delta_M$ is the Laplace-Beltrami operator given by
\begin{equation*}
    \Delta_Mu = \frac{1}{\mathbf{g}(u)}\mathrm{div}\left(\mathbf{g}(u) (u_x^{\top}u_x)^{-1} u_x^{\top}\right)
    \end{equation*}
with Gramian determinant $\mathbf{g}(u)= \det\left(u_x^{\top} u_x\right)^{1/2}$; see \cite[Ch.3.1]{Br78},
\cite{Sinestrari2010}.
The Laplace-Beltrami operator is known to be equivariant with respect to 
$\mathrm{Diff}(M)$, so that equation \eqref{eq:mcf} is already in 
coordinate-free form. The resulting Freezing system \eqref{e37a}, \eqref{e37b}
then keeps the type  of a parabolic system with quasilinear elliptic principal part. 
\newline

\textbf{Acknowledgment:} This paper is
based in part on the thesis \cite{Roe19} which was
initiated and supervised  by WJB and YGK.  We have
incorporated some theoretical advances,  and MEK has reorganized the numerical approach and carried out extended
numerical experiments.
The work of YGK was supported in part by the US National Science Foundation and the US Department of Energy.  YGK also wishes to  gratefully acknowledge decades-long discussions and collaboration along this conceptual direction with Prof. P. G. Kevrekidis of the University of Massachusetts, starting with their joint work on PDE focusing in the early 2000 years. This ongoing collaboration has touched, at times, on several aspects of the present work, and is currently being separately prepared for publication. 


  \bibliographystyle{abbrv}

\appendix
\renewcommand{\thesection}{Appendix \Alph{section}}
\renewcommand{\theequation}{\Alph{section}.\arabic{equation}}
\renewcommand{\thefigure}{\Alph{section}.\arabic{figure}}
\setcounter{equation}{0}
\setcounter{figure}{0}
\section{Comments on the analytical framework}
\label{sec2.3}

In this Appendix we indicate how 
the formal computations of Sections \ref{sec1} and \ref{sec2} can be put
into a mathematically rigorous framework. 
Further, we want to link up the Freezing system \eqref{e37a}, \eqref{e37b}
to the well-known approach for finite dimensional Lie groups;
see \cite{BOR14}, \cite{BO19}. A more detailed elaboration
is given in the thesis \cite{Roe19}.

The main  source for our undertaking is the monograph \cite{KM97} on global analysis.
The authors of \cite{KM97} work in the category of
{\it convenient vector spaces} which are locally convex topological
vector spaces satisfying a weak completeness condition, see \cite[2.14]{KM97}.
For example, Fr\'{e}chet spaces (cf. \cite[2.8.1]{J81}, \cite[p.48]{S71})
are special types of convenient vector spaces, see \cite[2.14]{KM97}.
In particular, our phase space 
$\mathcal{U}_d = C^{\infty}(M,\R^d)$ where $M$ is a finite dimensional manifold
becomes a  Fr\'{e}chet space when endowed with the topology of
uniform convergence of all derivatives on  compact subsets, see
\cite[52.9]{KM97}.


According to \cite[Theorem 43.1]{KM97} the group
$G=\mathrm{Diff}(M)$ of $C^{\infty}$-diffeo\-morphisms of the manifold $M$ is a 
Lie group, i.e. inversion and composition are smooth operations with respect
to the topology on the space of all smooth mappings $\mathfrak{C}^{\infty}(M,M)$,
see \cite[42.1]{KM97}.
The Lie algebra $\mathfrak{g}$ of $\mathrm{Diff}(M)$ is defined as the
tangent space $T_{\one}(G)=\mathfrak{g}$ at unity $\one=\mathrm{id}_M$ and agrees with the  convenient
vector space of all smooth vector fields on $M$ with compact support, denoted
by  $C_c^{\infty}(M \leftarrow TM)$ (see \cite[Theorem 43.1]{KM97}).
Since $M$ is compact in our setting, the space $C_c^{\infty}(M \leftarrow TM)$
agrees with the space $\mathcal{V}=C^{\infty}(M,TM)$ from \eqref{e2a}.
It is a well accepted fact that such a calculus of infinite-dimensional Lie groups
does not exist for manifolds of mappings with finite smoothness, for example in $C^k (M,M)$, see \cite[Introduction]{KM97}, \cite{Sch10}.

In the calculus above, the derivative of the multiplication from the left
\begin{equation*} \label{e18}
L_{\gamma}: G \rightarrow G, \quad L_{\gamma}(g) = \gamma \circ g, \; g \in G.
\end{equation*}
 at $g=\one$ turns out to be
\begin{equation} \label{e20}
d_gL_{\gamma}(\one)\mu= \partial_x\gamma \circ \mu  \quad \text{for} \quad 
\mu \in \mathfrak{g}=C^{\infty}(M,TM).
\end{equation}
Note that $\mu(x) \in T_xM$ for $x \in M$ and $\partial_x\gamma(x):
T_xM \to T_{\gamma(x)}M$, hence one has $\partial_x\gamma(x) \circ \mu(x) \in T_{\gamma(x)}M$ and the derivative $d_gL_{\gamma}(\one)$ maps the tangent space
$T_{\one}G=\mathfrak{g}$ into the tangent space $T_{\gamma}G=T_{L_{\gamma}(\one)}G$ as expected.

In a similar way,  the derivative $\partial_ga(g,u)$ of the action from
\eqref{e6} at  $g=\one$ is found to be
\begin{equation} \label{e22}
\partial_ga(\one,u) \mu= - \partial_xu \circ \mu, \quad \mu \in \mathfrak{g} = C^{\infty}(M,TM).
\end{equation}
The map $\partial_ga(\one,u):\mathfrak{g} \rightarrow C^{\infty}(M,\R^d)=\mathcal{U}_d$
is also called the {\it generator of the group action on $\mathcal{U}_d$}.
With the relations \eqref{e20} and \eqref{e22} to the abstract setting
we may write the Freezing system
\eqref{e37a}, \eqref{e37b} as follows:

\begin{equation*}
\begin{aligned} 
  \partial_tv(t,\cdot)= & F(v(t,\cdot))- \partial_ga(\one,v(t))\mu(t) ,\quad v(0,\cdot)=u_0(\cdot), \\
\partial_t g(t,\cdot) = & d_gL_{g(t)}(\one) \mu(t), \quad g(0,\cdot)= \one. 
\end{aligned}
\end{equation*}
This  is precisely the form of the Freezing equations derived for finite
dimensional Lie groups in  \cite{BOR14}, \cite{BT04}.  Hence, using this calculus, our
new system \eqref{e37a}, \eqref{e37b} appears as a direct generalization of the
finite dimensional case.  

\end{document}